\documentclass[12pt,a4paper]{amsart}
\usepackage{amsmath, amsthm, amssymb, enumerate, mathtools}

\usepackage[breaklinks=true,colorlinks=true,linkcolor=blue,citecolor=red,urlcolor=green]{hyperref}

\usepackage[margin=1in,footskip=0.25in]{geometry}

\usepackage{color}

\usepackage[utf8]{inputenc}
\usepackage{eepic,epic}
\usepackage{esint}
\usepackage{multirow}
\usepackage{array}
\usepackage[utf8]{inputenc}
\usepackage{graphicx}
\usepackage{caption}
\usepackage{booktabs}
\usepackage{siunitx}
\usepackage{makecell}
\usepackage{mathrsfs}
\usepackage{enumitem}

\usepackage{theoremref}

\def\to{\longrightarrow}
\def\mapsto{\longmapsto}

\def\cwedge{\bigcirc\kern-1.07em\wedge\ }

\newcommand \ve{\varepsilon}
\newcommand \px{p}
\newcommand \qx{q}
\newcommand \gq{\gamma}
\newcommand \ut{\vt}
\def\vt{{\tau\hskip-4.55pt\iota\hskip.6pt}} 
\def\hz{\mathcal{H}}
\def\lz{\mathcal{L}}
\def\vz{\mathcal{V}}

\newcommand{\txm}{T\hskip-3pt_x\w\hn M}
\newcommand{\df}{d\hskip-.8ptf}

\newcommand{\tr}{\mathrm{tr}}

\theoremstyle{plain}

\newtheorem{thm}{Theorem}[section]
\newtheorem{lem}[thm]{Lemma}
\newcommand{\bal}{\begin{aligned}}
\newcommand{\eal}{\end{aligned}}

\usepackage{float}

\theoremstyle{definition}
\newtheorem*{defn*}{Definition}

\theoremstyle{remark}
\newtheorem{rem}[thm]{Remark}
\newtheorem*{pf}{Proof}

\numberwithin{equation}{section}

\newcommand{\w}{^{\phantom i}}

\newcommand{\az}{a}
\newcommand{\bz}{b\hs}

\newcommand{\hyp}{\hskip.5pt\vbox
{\hbox{\vrule width2.5ptheight0.5ptdepth0pt}\vskip2pt}\hskip.5pt}
\newcommand{\hs}{\hskip.7pt}
\newcommand{\hh}{\hskip.4pt}
\newcommand{\hn}{\hskip-.4pt}
\newcommand{\nh}{\hskip-.7pt}
\newcommand{\nnh}{\hskip-1pt}
\newcommand{\ric}{\mathrm{r}}
\newcommand{\sch}{\mathrm{h}}
\newcommand{\sca}{\mathrm{s}}
\newcommand{\ein}{\mathrm{e}}
\newcommand{\trc}{\mathrm{trc}}
\newcommand{\kf}{\omega}
\newcommand{\rf}{\rho}
\newcommand{\ef}{\eta}
\newcommand{\bbR}{\mathrm{I\!R}}
\newcommand{\bbC}{{\mathchoice {\setbox0=\hbox{$\displaystyle\mathrm{C}$}
\hbox{\hbox to0pt{\kern0.4\wd0\vrule height0.9\ht0\hss}\box0}} 
{\setbox0=\hbox{$\textstyle\mathrm{C}$}\hbox{\hbox 
to0pt{\kern0.4\wd0\vrule height0.9\ht0\hss}\box0}} 
{\setbox0=\hbox{$\scriptstyle\mathrm{C}$}\hbox{\hbox 
to0pt{\kern0.4\wd0\vrule height0.9\ht0\hss}\box0}} 
{\setbox0=\hbox{$\scriptscriptstyle\mathrm{C}$}\hbox{\hbox 
to0pt{\kern0.4\wd0\vrule height0.9\ht0\hss}\box0}}}} 

\def\bbCP{\bbC\mathrm{P}}
\def\bbRP{\bbR\mathrm{P}}

\begin{document}	
\title[On weakly Einstein K\"ah\-ler surfaces] 
	{On weakly Einstein K\"ah\-ler surfaces}

\author[A.\,Derdzinski, Y\nnh.\,Euh, S.\,Kim, J,\,H.\,Park]{Andrzej
Derdzinski$^{1}$, Yunhee Euh$^{2}$, Sinhwi Kim$^{2}$, JeongHyeong Park$^{2}$}
    \address{$^{1}$Department of Mathematics, The Ohio State University, Columbus, OH
43210, USA}
	\address{$^{2}$Department of Mathematics, Sungkyunkwan University, Suwon, 16419, Korea}

\email{andrzej@math.ohio-state.edu, prettyfish@skku.edu, kimsinhwi@skku.edu,\\
parkj@skku.edu}


\begin{abstract}
Riemannian four-manifolds in which the triple contraction of the curvature
tensor against itself yields a functional multiple of the metric are called 
\emph{weakly Ein\-stein}. We focus on weakly Ein\-stein K\"ah\-ler surfaces.
We provide several conditions characterizing those K\"ah\-ler surfaces which 
are weakly Ein\-stein, classify weakly Ein\-stein K\"ah\-ler surfaces having
some specific additional properties, and construct new examples.
\end{abstract}
	
\maketitle

\section{Introduction}\label{in}
Following Euh, Park and Sekigawa \cite{EPS13}, we say that a Riemannian
four-man\-i\-fold is {\it weakly Ein\-stein\/} when the three-in\-dex
contraction of
its curvature tensor against itself equals a functional multiple of the metric.
This is the case -- in dimension four only -- for all Ein\-stein manifolds.
See formulae (\ref{def}) and (\ref{trf}) in Sect.\,\ref{ac}.

The same requirement in dimensions $\,n>4$, coupled with the Ein\-stein 
condition, defines what one calls {\it super-Ein\-stein manifolds}, and then
the above `functional multiple' must -- only if $\,n>4\,$ -- be 
a {\it constant multiple\/} \cite[p.\,165]{besse},
\cite[p.\,358]{GW}, \cite[Lemma 3.3]{BV}.

The present paper deals with weakly Ein\-stein K\"ah\-ler surfaces. 
Our Theorem~\ref{equiv} provides five conditions
equivalent to the weakly Ein\-stein property of a given K\"ah\-ler surface,
which allows us to quickly conclude whether a specific K\"ah\-ler surface is
weakly Ein\-stein.

One calls a function $\,\vt\,$ on a Riemannian manifold {\it
trans\-nor\-mal\/} if the integral curves of its
gradient $\,v\,$ are re\-pa\-ram\-e\-trized geodesics. This amounts to
requiring that, 
locally, at points where $\,v\ne0$, the norm-squar\-ed $\,Q\,$ of $\,v\,$ be a 
function of $\,\vt$. When the last condition holds for both $\,Q\,$ and the
La\-plac\-i\-an of $\,\vt$, one refers to $\,\vt\,$ as {\it
  iso\-par\-a\-met\-ric}. See \cite{bolton,miyaoka,wang}. Our next four 
theorems (three local, one assuming compactness) classify 
weakly Ein\-stein K\"ah\-ler surfaces having specific additional properties.
\begin{thm}\label{lchmg}Up to local isometries, Riemannian products of two
real surfaces of opposite nonzero constant Gauss\-i\-an curvatures are the
only locally homogeneous non-Ein\-stein weakly Ein\-stein K\"ah\-ler surfaces.
\end{thm}
\begin{thm}\label{sdlsy}All self-dual weakly Ein\-stein K\"ah\-ler surfaces
are locally symmetric, and hence locally homothetic to a standard complex 
projective or hyperbolic space, or to a product of two real surfaces of
opposite constant Gaussian curvatures.
\end{thm}
\begin{thm}\label{triso}On a weakly Ein\-stein K\"ah\-ler
surface, every trans\-nor\-mal function with a hol\-o\-mor\-phic gradient is
necessarily also iso\-par\-a\-met\-ric.
\end{thm}
\begin{thm}\label{trhol}Up to bi\-hol\-o\-mor\-phic isometries, the only
compact non-Ein\-stein weakly Ein\-stein K\"ah\-ler surfaces admitting 
nonconstant trans\-nor\-mal functions with hol\-o\-mor\-phic gradients are 
certain compact isometric quotients of the Riemannian product of a sphere and
a hyperbolic plane with opposite constant Gauss\-i\-an curvatures.
\end{thm}
The word `certain' in Theorem~\ref{trhol} accounts for the requirement that
both factor distributions be orient\-able (so as to make the quotient
K\"ah\-ler, and not just locally K\"ah\-ler). The isometric quotients in
question are precisely all the non\-flat con\-for\-mal\-ly flat compact
K\"ah\-ler surfaces. See Remark~\ref{cnffl}.

In the K\"ah\-ler-Ein\-stein case the conclusion of Theorem~\ref{triso} holds 
for well-known and general reasons. See the second sentence of
Remark~\ref{calab}.

We do not know whether there exist compact non-Ein\-stein weakly Ein\-stein
K\"ah\-ler surfaces other than the compact isometric quotients of the
Riemannian products of spheres and hyperbolic planes with opposite
Gauss\-i\-an curvatures, mentioned in Theorem~\ref{trhol}.

What we can easily conclude, using Theorem~\ref{trhol}, is that examples of
this kind -- if they exist -- cannot have a nonzero Euler characteristic and
at the same time admit groups of isometries with an infinite center and with
principal orbits of dimension three. Thus, a
$\,\mathrm{U}\hh(2)$-in\-var\-i\-ant K\"ah\-ler metric on the 
one-point and two-point blow-ups of $\,\bbCP^2$ is never weakly Ein\-stein.

Theorem~\ref{triso}, although not {\it per se\/} a classification result,
provides a crucial step in the proof of our Theorem~\ref{weskr}, which 
explicitly describes the local structure of all 
non-Ein\-stein weakly Ein\-stein K\"ah\-ler surfaces with 
nonconstant trans\-nor\-mal functions having hol\-o\-mor\-phic gradients. We
then use Theorem~\ref{weskr} to prove Theorem~\ref{trhol}.

As an added bonus, Theorem~\ref{weskr} leads to new examples of (noncompact) 
non-Ein\-stein weakly Ein\-stein K\"ah\-ler surfaces, presented
in Sect.\,~\ref{ne}.

For compact Riemannian $\,n$-manifolds with  parallel Ric\-ci tensor,
$\,n\ge4$, the weakly Ein\-stein property (as stated above when $\,n=4$) 
is necessary and sufficient in order that the metric be a critical point of 
the functional associating with metrics of unit volume the $\,L\nh^2$
norm-squared of their curvature \cite[Corollary 4.72]{besse}. See also 
\cite{EPS15}. For $\,n=3$, non-Ein\-stein weakly Ein\-stein manifolds are
characterized by having a Ric\-ci tensor of rank one \cite{GHMV}. For
$\,n\geq 4$, only partial classification results exist such as
Arias-Marco and Kowalski's theorem \cite{AK} in the locally homogenous case
with $\,n=4$, mentioned in Sect.\,\ref{ph}. Several results in this direction
can be found in \cite[Sect,\,6.55--6.63]{besse}.
Con\-for\-mal\-ly flat weakly Ein\-stein Riemannian manifolds were
classified by Garc\'\i a-R\'\i o et al.\ \cite{GHMV}. For the
weakly Ein\-stein
condition in extrinsic geometry, see \cite{KNP} and the references therein.

\section{Notation and preliminaries}\label{pr}
All manifolds, mappings, tensor fields and connections are assumed smooth.
Manifolds are by definition connected. 
Given a Riemannian manifold $\,(M\nh,g)$, we denote by
$\,\nabla\nnh,R,\ric,\ein\,$ and $\,\sca\,$ its Le\-vi-Ci\-vi\-ta connection,
curvature tensor, Ric\-ci tensor, Ein\-stein tensor and scalar curvature,
with the sign convention such that
$\,\ric(w,w')=\mathrm{tr}\hs[R(w,\,\cdot\,)w']\,$ for vector fields
$\,w,w'\nh$. Thus, $\,\sca=\tr\nh_{g}\w\hh\ric$. The symbol $\,\nabla\hn$ also
stands for the $\,g$-grad\-i\-ent.

In the underlying real space of a complex 
vector space $\,V\hs$ of positive finite dimension $\,m\,$ we use the natural 
orientation such that $\,e_1,ie_1,\dots,e_m,ie_m$ is a positive real basis 
for any complex basis $\,e_1,\dots,e_m$ of $\,V\nh$. (The automorphism group 
$\,\hs\text{\rm GL}\hs(V)\,\approx\,\text{\rm GL}\hs(m,\bbC)\,$ is connected,
since every automorphism has, in some basis, a triangular matrix, which 
is joined to $\,\hs\text{\rm Id}\hs\,$ by an obvious curve 
of nonsingular triangular matrices.) Thus,
\begin{equation}\label{nor}
\mathrm{every\ almost}\hyp\mathrm{complex\ manifold\ is\ canonically\
oriented.}
\end{equation}

\begin{rem}\label{seeqz}
A Riemannian manifold $\,(M\nh,g)\,$ with $\,\sca\hh\ein=0\,$ that is,
one in which $\,\sca\hh\ric\,$ is a functional multiple of $\,g$, necessarily
has $\,\sca=0\,$ identically, or $\,\ein=0\,$ everywhere.
(This needs to be justified as we are not assuming 
real-an\-a\-lyt\-ic\-i\-ty.) First, we may assume that $\,\dim M>2$. (For
surfaces, $\,\sca\hh\ein=\ein=0$.) Now
$\,\sca\hskip1.7ptd\hh\sca=0$, since a point with $\,\sca\ne0\,$ has a
connected neighborhood $\,\,U\,$ on which $\,\ein=0$, so that, by Schur's
lemma, $\,\sca\,$ is constant and $\,d\hh\sca=0\,$ on $\,\,U\nh$. The
ensuing constancy of $\,\sca^2$ gives $\,\sca=0\,$ identically or
$\,\sca\ne0\,$ everywhere.
\end{rem}
\begin{rem}\label{posdr}If $\,d<c\,$ and a differentiable function $\,\psi\,$
of the variable $\,\alpha\in[d,c)\,$ has a positive derivative 
$\,\psi\hs'(\alpha)\,$ at every $\,\alpha\,$ at which
$\,\psi(\alpha)\ge\lambda$,
for some given $\,\lambda\in\bbR$, then either $\,\psi<\lambda\,$ everywhere
in $\,[d,c)\,$ or, for the least $\,\alpha_1\w\in[d,c)\,$ at which
$\,\psi(\alpha_1\w)\ge\lambda$, the restriction of $\,\psi\,$ to 
$\,(\alpha_1\w,c)\,$ is greater than $\,\lambda\,$ and strictly increasing;
consequently, $\,\psi\,$ has at $\,c\,$ a limit lying in
$\,(\lambda,\infty]$.

In fact, our claim will follow immediately once we verify that
$\,\psi\ge\lambda\,$ on $\,(\alpha_1\w,c)$. Assuming this not to be the case, 
we let
$\,\alpha_2\w$ be the infimum of those $\,\alpha\in(\alpha_1\w,c)\,$ at which
$\,\psi(\alpha)<\lambda$. Thus, $\,\alpha_2\w>\alpha_1\w$, since at
$\,\alpha_1\w$ one has $\,\psi(\alpha_1\w)\ge\lambda\,$ and
$\,\psi\hs'(\alpha_1\w)>0$. It follows now that $\,\psi(\alpha_2\w)=\lambda$,
and hence $\,\psi\hs'(\alpha_2\w)>0$, giving $\,\psi(\alpha)>\lambda\,$ for
all $\,\alpha>\alpha_2\w$ close to $\,\alpha_2\w$, contrary to the definition
of $\,\alpha_2\w$.
\end{rem}

\section{Proof of Theorem~\ref{lchmg}}\label{ph}
Recall from the Introduction that a Riemannian four-man\-i\-fold is said to be 
weakly Ein\-stein when the three-in\-dex contraction of its curvature tensor
$\,R\,$ against itself equals some function $\,\phi\,$ times the metric
$\,g\,$ (in coordinates: $\,R_{ikpq}\w R_j\w{}^{kpq}\nnh=\phi g_{ij}\w$).

Arias-Marco and Kowalski \cite{AK} showed that a non-Ein\-stein locally
homogeneous weakly Ein\-stein four-manifold must be locally isometric
either to a Riemannian product of 
surfaces with opposite nonzero constant Gauss\-i\-an curvatures, or to what
they call an EPS space: one of the examples constructed by 
Euh, Park, and Sekigawa \cite[Example 3.7]{EPS14}.

This proves Theorem~\ref{lchmg}, since an EPS space $\,(M\nh,g)\,$ is not 
K\"ah\-ler.

We verify the last claim as follows. In \cite[Example 3.7]{EPS14},
$\,M\,$ is a Lie group with left-in\-var\-i\-ant 
$\,g$-or\-tho\-nor\-mal vector fields $\,e_1\w,\dots,e_4\w$ satisfying some  
Lie-brack\-et relations 
that involve constants $\,a\ne0\,$ and $\,b$. According to 
\cite[formula (3.14)]{EPS14}, the only nonzero components 
of the curvature tensor $\,R\,$ in this frame are those algebraically related 
to $\,R_{1212}\w=R_{1313}\w=R_{1414}\w=R_{3434}=-a^2$ and  
$\,R_{2323}\w=R_{2424}\w=a^2$. (The sign of $\,R\,$ in 
\cite{EPS14} is the opposite of ours.) Consequently, the frame
di\-ag\-o\-nal\-iz\-es the Ric\-ci tensor $\,\ric$, with the eigen\-val\-ues 
$\,\ric_{11}\w=-3a^2\nh$, $\,\ric_{22}\w=a^2$ and
$\,\ric_{33}\w=\ric_{44}\w=-a^2\nh$. 

Thus, $\,g\,$ is not a K\"ah\-ler metric: if it were, the complex structure,
leaving $\,\ric\,$ invariant, would cause $\,\ric\,$ to have two double
eigen\-val\-ues, contrary to the above formulae.

\section{Algebraic curvature tensors}\label{ac}
For an algebraic curvature tensor $\,R\,$ in 
a Euclidean space $\,\mathcal{T}\hs$ of dimension $\,n\ge4$, denoting by
$\,g\,$ the inner product, by $\,W\nnh\nh,\ric,\ein\,$ the Weyl,
Ric\-ci and Ein\-stein tensors of $\,R$, and by $\,\sca\,$ its scalar
curvature, with
$\,\ric_{ij}\w=g^{pq}R_{ipjq}\w$, $\,\sca=g^{pq}\ric_{pq}\w$ and
$\,\ein=\ric-\sca g/n$, we have
\begin{equation}\label{wey}
\begin{array}{rl}
\mathrm{a)}&W_{\!ijpq}\w=R_{ijpq}\w
-\displaystyle{\frac1{n-2}}(g_{ip}\w\ric_{jq}\w
+g_{jq}\w\ric_{ip}\w-g_{jp}\w\ric_{iq}\w
-g_{iq}\w\ric_{jp}\w)\\
&\hskip69.5pt
+\,\displaystyle{\frac{\sca}{(n-1)(n-2)}}
(g_{ip}\w g_{jq}\w-
g_{jp}\w g_{iq}\w).\phantom{1^{1^1}}\\
\mathrm{b)}&W_{\!ijpq}\w=R_{ijpq}\w
-\displaystyle{\frac1{n-2}}(g_{ip}\w\sch_{jq}\w
+g_{jq}\w\sch_{ip}\w-g_{jp}\w\sch_{iq}\w
-g_{iq}\w\sch_{jp}\w),
\end{array}
\end{equation}
where the more concise version involves the Schou\-ten tensor
$\,\sch=\ric-\sca g/(2n-2)$, 
We use here components relative to any basis of $\,\mathcal{T}\nh$, with
index raising and lowering via $\,g$, and summation over
repeated indices. See \cite[Sect.\,1.108]{besse}.

The {\it triple contraction\/} of {\it any such\/} $\,R\,$ is the symmetric
$\,2$-ten\-sor $\,\trc\hs R\,$ given by
\begin{equation}\label{trc}
[\trc\hs R]_{ij}\w\,=\,\,R_{ikpq}\w R_j\w{}^{kpq}\nh.
\end{equation}
Using (\ref{wey}-b) and, respectively, (\ref{wey}-a), one easily verifies that
\begin{equation}\label{trw}
\begin{array}{rl}
\mathrm{i)}&(n-2)^2(\trc\hs R\,-\,\trc\,W)\,
+\,2[2\sca\hh\ric\,-\,2(n-2)R\hh\ric\,-\,n\ric^2\hh]\,\mathrm{\ is\ a\
multiple\ of\ }\,g,\\
\mathrm{ii)}&(n-2)(R\hh\ric\,-\,W\nnh\ric)\,+\,2\ric^2\hh\,
-\,n\sca\hh\ric/(n-1)\,\mathrm{\ is\ a\ multiple\ of\ }\,g.
\end{array}
\end{equation}
Here $\,b=\ric^2$ has the components $\,b_{ij}\w=\ric_j^k\ric_{ik}\w$,
and an algebraic curvature tensor $\,R\,$ acts on arbitrary 
$\,(0,2)\,$ tensors $\,b\,$ by $\,[Rb]_{ij}\w=R_{ipjq}\w b\hs^{pq}\nh$, which
preserves (skew)sym\-metry of $\,b$, and --  due to the Bianchi identity -- 
becomes $\,2[R\alpha]_{ij}\w=R_{ijpq}\w\alpha^{pq}$ when $\,b=\alpha\,$ 
is skew-sym\-met\-ric. 
We will repeatedly refer to tensors (or, tensor fields on a manifold with a
metric $\,g$) as multiples (or, functional multiples) of $\,g$, without the
need to specify the factor, since it is trivially found by contraction.

For instance, given a Riemannian manifold $\,(M\nh,g)\,$ with the
curvature tensor $\,R$,
\begin{equation}\label{def}
(M\nh,g)\,\mathrm{\ is\ called\ }\mathit{weakly\ Einstein}\mathrm{\ when
\ }\,\trc\hs R\,\mathrm{\ \ is\ a\ functional\ multiple\ of\ }\,g.
\end{equation}
It is well known -- see, e.g., \cite[p.\,413]{cm83} -- that $\,\trc\,W\,$ 
is a multiple of $\,g\,$ when $\,n=4$. Thus, (\ref{trw}-i) with $\,n=4\,$
implies that, as pointed out by Euh, Park and Sekigawa \cite{EPS13},
\begin{equation}\label{trf}
\begin{array}{l}
\mathrm{in\ dimension\ four,\ \ \ 
}\,\trc\hs R-2R\hh\ric+\sca\hh\ric-2\ric^2\hs\mathrm{\ \ \ is\
a\ multiple\ of\ \ }\,g\mathrm{,\ so}\\
\mathrm{that,\hskip1.5ptfor\ Riemannian\
}\hs4\hyp\mathrm{manifolds,\hskip1.5ptEinstein\
implies\ weakly\ Einstein.}
\end{array}
\end{equation}
Setting $\,n=4\,$ in (\ref{trw}-ii) and adding the resulting expression to
the one in (\ref{trf}), 
\begin{equation}\label{trm}
\mathrm{we\ see\ that\ }\,\trc\hs R-2\hh W\nnh\ric-\sca\hh\ric/3\,\mathrm{\
is\ a\ multiple\ of\ }\,g\,\mathrm{\ when\ }\,n=4,
\end{equation}
which easily yields \cite[Theorem 2(i)]{GHMV}: 
a four-di\-men\-sion\-al non-Ein\-stein con\-for\-mal\-ly flat Riemannian
manifold is weakly Ein\-stein if and only if its scalar curvature vanishes.
(This involves a technical detail, discussed in Remark~\ref{seeqz}.) As the
contractions of $\,W\nnh\nh$ vanish, which makes $\,W\nnh\ric\,$ obviously
equal to $\,W\nnh\nh\ein$, (\ref{trm}) implies that, in dimension four,
\begin{equation}\label{iff}
(M\nh,g)\,\mathrm{\ is\ weakly\ Ein\-stein\ if\ and\ only\ if\
}\,\,6\hh W\nnh\nh\ein=-\sca\hh\ein\hh,
\end{equation}
$\ein=\ric-\sca g/4\,$ being again the Ein\-stein tensor. 
The additional assumption
\begin{equation}\label{add}
\begin{array}{l}
\mathrm{that\ }\,\,\ein\,\,\mathrm{\ has\ the\ spectrum\ of\ the\ form\ \ 
}\,\hs(a,\hs a,-a,-a),\\
\mathrm{which\ is\ obviously\ satisfied\ by\
K}\ddot{\mathrm a}\mathrm{h\-ler\ surface\ metrics,}
\end{array}
\end{equation}
leads to the following consequence:
\begin{equation}\label{trq}
\sca\hh\ric-2\ric^2\hs\mathrm{\ is\ a\ multiple\ of\ }\,g.    
\end{equation}
In fact, by (\ref{add}), $\,\ein^2$ is a multiple of $\,g$, while
$\,2\ein^2\nh=2\ric^2\nh-\sca\hh\ric+\sca^2g/8$.

Therefore, from (\ref{trf}) -- (\ref{trq}) it trivially follows that
\begin{equation}\label{rrr}
\trc\hs R-2R\hh\ric\,\mathrm{\ is\ a\ multiple\ of\ }\,g\,\mathrm{\ if\ }\,n
=4\,\mathrm{\ and\ (\ref{add})\ holds.}
\end{equation}

\section{K\"ah\-ler surfaces}\label{ks}
A twice-co\-var\-i\-ant tensor field $\,\az\,$ on 
an al\-most-com\-plex manifold $\,(M\nh,J)\,$ gives rise to two more such 
tensor fields,
\begin{equation}\label{aje}
\,\az J=\az(J\hs\cdot\,,\,\cdot\,)\,\mathrm{\ \  and\ \ }\,
J\az=-\az(\,\cdot\,,\,J\hs\cdot\,),
\end{equation}
as well as the {\it commutator} 
$\,[\az,J]=-[J,\az]=\az J-J\az$. Obviously,
\begin{equation}\label{inv}
\az\,\mathrm{\ is\ }\,J\hyp\mathrm{invariant,\ that\ is,\
}\az(J\hs\cdot,J\hs\cdot\,)=\az\mathrm{,\ if\ and\ only\ if\ }\,\az J=J\az.
\end{equation}
The tensor field $\,\az\,$ is said to be {\it 
Her\-mit\-i\-an\/} if it is symmetric at every point and 
$\,J$-in\-var\-i\-ant. Clearly, if $\,\az\,$ is
Her\-mit\-i\-an, $\,\az J\,$ is a $\,2$-form and $\,(\az J)J=J(J\az)=-\hs\az$.
By a {\it Her\-mit\-i\-an metric\/} on an al\-most-com\-plex manifold
$\,(M\nh,J)\,$
we mean a Riemannian metric $\,g\,$ on $\,M\,$ which is a Her\-mit\-i\-an
tensor ($gJ=J\nh g$), which amounts to skew-ad\-joint\-ness of
$\,J\,$ at every point or -- equivalently -- the 
requirement that $\,J\,$ act in every tangent space as a linear isometry. 
If $\,g\,$ is Her\-mit\-i\-an, and one identifies bundle morphisms
$\,A,B:T\nh M\to T\nh M\,$ with $\,\az=g(A\hs\cdot\,,\,\cdot\,)\,$ and 
$\,\bz=g(B\hs\cdot\,,\,\cdot\,)$, then 
the operation $\,\az\mapsto\bz=J\az\,$ (or, 
$\,\az\mapsto\bz=\az J$), defined by (\ref{aje}) for 
twice-co\-var\-i\-ant tensor fields $\,\az$,
\begin{equation}\label{coi}
\mathrm{corresponds\ to\ the\ ordinary\ composition\
}\,B=J\nh A\,\mathrm{\ or,\ }\,B=AJ.
\end{equation}
As an example, in a K\"ah\-ler manifold $\,(M\nh,g,J)\,$ of real dimension
$\,n$, both $\,g\,$ and the Ric\-ci tensor $\,\ric\,$ are Her\-mit\-i\-an, and
hence so is the Ein\-stein (trace\-less Ric\-ci) tensor $\,\ein=\ric-\sca g/n$,
giving rise to the K\"ahler, Ric\-ci and Ein\-stein $\,2$-forms
$\,\kf=gJ,\,\rf=\ric J\,$ and $\,\ef=\ein J$, with
\begin{equation}\label{kre}
\nabla\nh\kf=0,\qquad d\rf=0,\qquad\ef=\rf-\sca\hs\hh\kf/n.
\end{equation}
In any complex manifold $\,M\nh$, the operator  
$\,i\hskip1pt\partial\overline{\partial}\,$ sends every 
$\,C^\infty$ function $\,f:M\to\bbR$ to the exact $\,2$-form 
$\,i\hskip1pt\partial\overline{\partial}\nh f\,$ given by
\begin{equation}\label{idd}
2\hs i\hskip1pt\partial\overline{\partial}\nh f\,=\,-\hs d\hskip1pt[(\df)J]\hs.
\end{equation}
Here the $\,1$-form $\,(\df)J\,$ equals, at any point $\,x\in M$, the 
composite of $\,J_x\w:\txm\to \txm\,$ followed by 
$\,\df\nnh_x\w:\txm\to\bbR\hs$. For any tor\-sion\-free connection
$\,\nabla\hs$ on the com\-plex manifold $\,M\,$ such that $\,\nabla\nh J=0$,
(\ref{idd}) is easily seen to become
\begin{equation}\label{jdj}
2i\hskip1pt\partial\overline{\partial}\hskip-.9ptf\hs
=\bz J\,+\,J\bz,\quad\mathrm{where\,}\,\bz=\nabla\nh\df\nh.
\end{equation}
The Weyl tensor acts on $\,2$-forms $\,\alpha\,$ in a Riemannian 
manifold of dimension $\,n\ge4$, with the scalar curvature $\,\sca$, via 
the Weitzen\-b\"ock formula \cite[p.\,409]{cm83}, \cite[p.\,458]{hdg00}, 
immediate from the Ric\-ci identity:
\begin{equation}\label{wei}
W\nnh\alpha=\displaystyle{\frac12}\left[\delta(\nabla\nh\alpha-d\alpha)
-d\hs\delta\alpha)\right]
+\displaystyle{\frac{n-4}{2(n-2)}}\,\{\ric,\alpha\}
+\displaystyle{\frac{\sca}{(n-1)(n-2)}}\hs\alpha,
\end{equation}
$\{\,,\}\,$ being the anticommutator, $\,\ric\,$ the Ric\-ci tensor. In local
coordinates,
\begin{equation}\label{lco}
\begin{array}{l}
W_{\!ijpq}\w\alpha^{pq}\nh=-\,\alpha_{pi,\hh j}\w{}^p-\alpha_{jp,\hh i}\w{}^p
-\alpha_{pj,}\w{}^p{}_i\w+\alpha_{pi,}\w{}^p{}_j\w\phantom{\displaystyle{\frac12}}\\
\phantom{W_{\!ijpq}\w\alpha^{pq}\nnh}
+\,\displaystyle{\frac{n-4}{n-2}}(\ric_j^{\hs p}\alpha_{ip}\w
+\ric_i^{\hs p}\alpha_{pj}\w)
+\displaystyle{\frac{2\hskip.8pt\sca}{(n-1)(n-2)}}\hs\alpha_{ij}\w.
\end{array}
\end{equation}
When $\,n=4$, (\ref{kre}) and (\ref{lco}) yield
\begin{equation}\label{twr}
\mathrm{a)}\hskip6ptW\nnh\kf\,=\,\displaystyle{\frac{\sca}6}\hs\kf,\qquad
\mathrm{b)}\hskip6pt2\hh W\nnh\nh\ef\,
=\,\Delta\rf\,
-\,i\hskip1pt\partial\overline{\partial}\hskip1.3pt\sca\,
+\,\displaystyle{\frac{\sca}3}\hs\ef.
\end{equation}
Here $\,\Delta\,$ sends a twice-co\-var\-i\-ant tensor field $\,\az\,$ 
to $\,\Delta\az\,$ given by $\,[\Delta\az]_{pq}\w=\az_{pq,k}\w{}^k\nh$. 
\begin{thm}\label{equiv}
For a K\"ah\-ler surface\/ $\,(M\nh,g,J)\,$ with the Weyl, Ric\-ci, Ein\-stein
tensors\/ $\,W\nnh\nh,\ric,\ein$, Ric\-ci and Ein\-stein forms\/ $\,\rf,\ef$,
and the scalar curvature $\,\sca$, the following six conditions are mutually
equivalent.
\begin{itemize}
\item[\rm{(a)}]$(M\nh,g,J)\,$ is weakly Ein\-stein.
\item[\rm{(b)}]$R\hh\ric\,$ equals a functional multiple of $\,g$.
\item[\rm{(c)}]$6\hh W\nnh\nh\ein\,=\,-\sca\hh\ein$.
\item[\rm{(d)}]$3W\nnh\nh\ef\,=\,-\sca\hh\ef$.
\item[\rm{(e)}]$\Delta\rf\,
-\,i\hskip1pt\partial\overline{\partial}\hskip1.3pt\sca\,
=\,-\sca\hs\ef$.
\item[\rm{(f)}]$2\hh\Delta\ric\,-\,\nabla\nh d\sca\,
+\,J[\nabla\nh d\sca]J\,=\,-2\sca\hs\ein$.
\end{itemize}
\end{thm}
We prove Theorem~\ref{equiv} in the next section.
\begin{rem}\label{justf}
Since $\,R\hh\ric\,$ behaves ``multiplicatively'' under Riemannian products, 
and equals $\,K^2\nh g\,$ for a surface metric with Gauss\-i\-an curvature
$\,K$, a Riemannian product of two surfaces, being locally K\"ah\-ler, is --
according to (a), (b) above -- weakly Ein\-stein if and only if the two
Gauss\-i\-an curvatures are constant and equal or mutually opposite. The
latter case constitutes the example found by Euh, Park, and Sekigawa
\cite{EPS14}.
\end{rem}
\begin{rem}\label{algct}
The equivalence of {\rm(a)} through {\rm(c)} above really amounts to 
a statement about algebraic curvature tensors, once we 
replace the K\"ah\-ler condition by {\rm(\ref{add})} and {\rm(a)} by
$\,\trc\hs R\,$ {\rm is a multiple of} $\,g$. We could include {\rm(d)} here 
as well if, instead of just {\rm(\ref{add})}, we invoked one further
property of K\"ah\-ler-type algebraic curvature tensors, namely,
{\rm(\ref{twr}-a)}. For more details, see Section\/~{\rm\ref{pt}}.
\end{rem}
\begin{rem}\label{lcsym}
As shown by Car\-tan \cite{cartan}, see also \cite[Theorem\,14.7]{hdg00}, 
lo\-cal-ho\-moth\-e\-ty types of locally symmetric Riemannian
four-man\-i\-folds form seven disjoint classes: {\it flat, spherical,
hyperbolic, complex projective, complex hyperbolic\/} and, finally,
{\it non\-flat\/ $\,2+2\,$ and\/ $\,1+3$ Riemannian products} (with
Ein\-stein factors). Of these, the first five -- due to their being
Ein\-stein, cf. (\ref{trf}) -- are weakly Ein\-stein. The sixth one is 
weakly Ein\-stein when the two Gauss\-i\-an curvatures 
are equal (the Ein\-stein case) or opposite (Remark~\ref{justf}).
Remark~\ref{justf} also shows that the remaining products listed
above are not weakly Ein\-stein, the $\,1+3\,$ case, with $\,W\nnh\nh=0\,$ 
and $\,\sca\ne0$, being immediate from the sentence following (\ref{trm}).
\end{rem}
\begin{rem}\label{cnffl}As shown by Tanno \cite{tanno}, up to local
isometries, the only con\-for\-mal\-ly flat K\"ah\-ler surfaces are the 
Riemannian products of real surfaces with opposite constant Gauss\-i\-an
curvatures.
\end{rem}
\begin{rem}\label{calab}
As observed by Ca\-la\-bi \cite{calabi}, cf.\ also
\cite[Sect.\,2.140]{besse}, on any K\"ah\-ler manifold, for 
a function $\,\vt\,$ with a real hol\-o\-mor\-phic gradient 
$\,v=\nabla\nh\vt$, one has $\,2\hh\ric(v,\,\cdot\,)=-d\Delta\vt$. In the
K\"ah\-ler-Ein\-stein case this implies that, locally, at points where
$\,v\ne0$, the La\-plac\-i\-an of $\,\vt\,$ is a function of $\,\vt$.
\end{rem}

\section{Proof of Theorem~\ref{equiv}}\label{pt}
That (a) is equivalent to both (b) and (c), is immediate: the former from 
(\ref{def}) and (\ref{rrr}), the latter due to (\ref{iff}). 
We now proceed to show that (c) holds if and only if (d) does, still using a
purely algebraic argument (cf.\ Remark~\ref{algct}), and adopting the notation
of the lines preceding (\ref{wey}), with $\,\mathcal{T}\hs$ standing for the
tangent space of $\,M\,$ at a given point. The inner product $\,g\,$ provides
the identifications
\begin{equation}\label{ide}
\mathrm{a)}\hskip6pt\mathcal{T}\hn=\mathcal{T}\hh^*\nh,\qquad 
\mathrm{b)}\hskip6pt\mathcal{T}^{\wedge2}\hn=\hs\mathfrak{so}\hs(\mathcal{T})
=[\mathcal{T}\hh^*]^{\wedge2}\nnh,
\end{equation}
so that $\,u\wedge v$, for
$\,u,v\in\mathcal{T}\nh$, becomes the endomorphism of $\,\mathcal{T}\hs$
given by
\begin{equation}\label{uvw}
w\,\mapsto\,(u\wedge v)w\,=\,\langle u,w\rangle v\,-\,\langle v,w\rangle u
\end{equation}
We assume that $\,\ein\ne0$, and choose in $\,\mathcal{T}\hs$ an
orthonormal basis of the form $\,u,Ju,v,Jv$, di\-ag\-o\-nal\-iz\-ing
$\,\ein\,$ with some eigen\-val\-ues $\,(a,a,-a,-a)$. Thus, by (\ref{uvw}),
$\,J=u\wedge Ju+v\wedge Jv$, and so, for the basis $\,\xi^1,\dots,\xi^4$
of $\,\mathcal{T}\hh^*$ dual to $\,u,Ju,v,Jv$,
\begin{equation}\label{jeu}
\kf=\xi^1\nnh\wedge\hs\xi^2+\hs\xi^3\nnh\wedge\hs\xi^4,\,\,\,
\ein=a(\xi^1\nnh\otimes\hs\xi^1+\hs\xi^2\nh\otimes\hs\xi^2
-\hs\xi^3\nh\otimes\hs\xi^3-\hs\xi^4\nh\otimes\hs\xi^4),\,\,\,
\ef=a(\xi^1\nnh\wedge\hs\xi^2-\hs\xi^3\nnh\wedge\hs\xi^4),
\end{equation}
with $\,a\ne0$. Here $\,\kf=g(J\hs\cdot\,,\,\cdot\,)\,$ and
$\,\ef=\ein(J\hs\cdot\,,\,\cdot\,)$, while, for 
$\,1$-forms $\,\xi,\zeta\in\mathcal{T}\hh^*\nnh$, we set
\begin{equation}\label{xtz}
[\xi\otimes\zeta]\hn_{pq}=\xi_p\w\zeta_q\w,\qquad[\xi\wedge\zeta]\hn_{pq}\w
=\xi_p\w\zeta_q\w-\xi_q\w\zeta_p\w.
\end{equation}
Due to (\ref{twr}-a) and (\ref{jeu}) with $\,a\ne0$, we have
$\,3W\nnh\nh\ef=-\sca\hh\ef\,$ if and only if 
\begin{equation}\label{wef}
12\hh W(\xi^1\nnh\wedge\hs\xi^2)=-\sca\,\xi^1\nnh\wedge\hs\xi^2\nh
+3\hh\sca\,\xi^3\nnh\wedge\hs\xi^4\nh,\quad
12\hh W(\xi^3\nh\wedge\hs\xi^4)=3\hh\sca\,\xi^1\nnh\wedge\hs\xi^2\nh
-\sca\,\xi^3\nnh\wedge\hs\xi^4\nh.
\end{equation}
As stated in the lines following (\ref{trw}), 
any algebraic curvature tensor $\,R\,$ acts on bivectors via
$\,2[R\alpha]_{ij}\w\hs=\,R_{ijpq}\w\alpha^{pq}\nh$, so that 
$\,R(w\wedge w')=R(w,w')$. Thus, (\ref{wef}) is nothing else than
\begin{equation}\label{amo}
\begin{array}{l}
W_{\!1212}\w=\hs W_{\!3434}\w=\hs-\sca/12\,\mathrm{\ \ and\ \
}\,W_{\!1234}\w=\hs\sca/4\mathrm{,\ \ while}\\
W_{\!12ij}\w\nh=W_{\!34ij}\w\nh=\,0\,\mathrm{\ unless\ }\,\{i,j\}\,\mathrm{\ 
is\ }\,\{1,2\}\,\mathrm{\ or\ }\,\{3,4\}.
\end{array}
\end{equation}
Since $\,[W\nnh\nh\ein]_{ij}\w=W_{\!ipjq}\w\ein\hs^{pq}\nh$, and
$\,W\nh\nnh g=0$, (\ref{jeu}) gives
$\,[W\nnh\nh\ein]_{11}\w\nh=[W\nnh\nh\ein]_{22}\w\nh=2a\hh W_{\!1212}\w$, 
$\,[W\nnh\nh\ein]_{12}\w\nh=[W\nnh\nh\ein]_{34}\w\nh=0$, 
$\,[W\nnh\nh\ein]_{13}\w\nh=2a\hh W_{\!1232}\w$, 
$\,[W\nnh\nh\ein]_{14}\w\nh=2a\hh W_{\!1242}\w$, 
$\,[W\nnh\nh\ein]_{23}\w\nh=2a\hh W_{\!1213}\w$,
$\,[W\nnh\nh\ein]_{24}\w\nh=2a\hh W_{\!1214}\w$, 
$\,[W\nnh\nh\ein]_{33}\w\nh=[W\nnh\nh\ein]_{44}\w\nh=-2a\hh W_{\!3434}\w$. As
$\,\ein_{11}\w\nh=\ein_{22}\w\nh=a=-\ein_{33}\w\nh=-\ein_{44}\w$ and
$\,\ein_{ij}\w\nh=0\,$ otherwise, this description of $\,W\nnh\nh\ein$,
combined with (\ref{amo}), proves that (c) is equivalent to (d). More
precisely, the above equalities amount to (\ref{amo}) except the formula
for $\,W_{\!1234}\w$ which, however, then follows as (\ref{twr}-a)
and (\ref{jeu}) yield $\,W_{\!1234}\w=-\nh W_{\!1212}\w+\hs\sca/6=\hs\sca/4$.

The equivalence of (d) and (e) is in turn immediate from (\ref{twr}-b).

Finally, (f) is -- due to (\ref{jdj}) -- precisely the result of applying
$\,J\,$ to (e).

\section{Proof of Theorem~\ref{sdlsy}}\label{po}
On any K\"ahler surface, with the orientation as in (\ref{nor}),
\begin{equation}\label{asd}
\mathrm{the\ Einstein\ form\ }\,\eta\,\mathrm{\ is\ always\
anti}\hyp\mathrm{self}\hyp\mathrm{dual.}
\end{equation}
This is clear since `anti-self-dual' is well known
\cite[Corollary\,37.3, Proposition\,37.5]{hdg00} to mean the same as
{\it orthogonal to\/ $\,\kf\,$ and commuting with\/} $\,\kf$, in the sense of 
the identification of $\,2$-forms with skew-ad\-joint endomorphisms provided
by (\ref{ide}-b), while $\,J$-in\-var\-i\-ance of $\,\ric\hs$ amounts to its
commuting with $\,J$.

Let a K\"ah\-ler surface $\,(M\nh,g)\,$ be weakly Ein\-stein and self-dual.

Theorem~\ref{equiv}(d) gives $\,3W\nnh\nh\ef=-\sca\hh\ef\,$ while,
from (\ref{asd}), $\,W\nnh\nh\ef=0\,$ due to the self-dual\-i\-ty assumption. 
Thus, $\,\sca\hh\ein=0\,$ everywhere. By Remark~\ref{seeqz}, one of
$\,\sca,\,\ein\,$ is identically zero, and in
either case $\,\sca\,$ is constant. 
Constancy of $\,\sca\,$ combined with self-dual\-i\-ty
implies local symmetry 
\cite[Proposition\ 9.3]{bourguignon}, \cite[Lemma\ 7]{cm83}.

The final clause of Theorem~\ref{sdlsy} 
is immediate from Remarks~\ref{lcsym} and~\ref{cnffl}.

\section{Proof of Theorem~\ref{triso}}\label{pm}
In \cite[Sect.\,5]{iumj12} one fixes a nonuple
$\,I\nh,a,\varSigma,h,\lz,\,(\,,\hs),\hz,\gamma,Q\,$ consisting of
\begin{itemize}
\item[\rm{(i)}] a nontrivial closed interval
$\,I\nh=[\vt\nh_{\mathrm{min}}\w,\vt\nh_{\mathrm{max}}\w]\,$ of the variable
$\,\vt$,
\item[\rm{(ii)}] a real number $\,a>0$,
\item[\rm{(iii)}] a compact K\"ah\-ler manifold $\,(\varSigma,h)\,$ of
complex dimension $\,1$,
\item[\rm{(iv)}] a function $\,Q:I\to\bbR\,$ equal to $\,0\,$ at the endpoints
of $\,I\hs$ and positive on its interior $\,I\nh^o\nnh$, such that 
$\,\dot Q(\vt\nh_{\mathrm{min}}\w)=2a=-\dot Q(\vt\nh_{\mathrm{max}}\w)$,
\item[\rm{(v)}] a mapping $\,\gamma:I\nh\to\bbRP^1\nh$, with
$\,I\nh\subseteq\bbR\subseteq\bbRP^1\nh$,
\item[\rm{(vi)}] a complex line bundle $\,\lz\,$ over $\,\varSigma\,$ with a
Her\-mit\-i\-an fibre metric $\,(\,,\hs)$,
\item[\rm{(vii)}]the horizontal distribution $\,\hz\hs$ of a connection in
$\,\lz\,$ making the fibre metric $\,(\,,\hs)$
parallel and having the curvature form
$\,-a(\vt\nnh_*\w\nh-\gamma)^{-\nh1}\omega^{(h)}\nnh$,
\end{itemize}
where $\,\vt\nnh_*\w\in I\hs$ is the midpoint, 
$\,\omega^{(h)}$ is the K\"ah\-ler form of $\,(\varSigma,h)\,$
and $\,(\,\,)\dot{\,}=\,d/d\vt$. One also fixes a dif\-feo\-mor\-phism 
$\,I\nh^o\nh\nh\ni\vt\mapsto r\in(0,\infty)\,$ such that
$\,\dot r=ar/Q$, and uses the symbol $\,r\hs$ both for an independent variable
ranging over $\,[\hh0,\infty)\,$ and for the norm function
$\,r:\lz\to[\hh0,\infty)\,$ of the fibre metric $\,(\,,\hs)$, so that our
fixed 
dif\-feo\-mor\-phism turns $\,\vt$, and hence $\,Q\,$ as well, into a
function on the total space $\,\lz$.

For the section $\,v\,$ of the vertical distribution $\,\vz\hs$ on $\,\lz\,$
which, restricted to each fibre, equals our $\,a\,$ times the radial
(identity) vector field, $\,d_v\w$ acts on functions of $\,\vt\,$ as
$\,Q\,d/d\vt$, and $\,v\,$ equals the $\,g$-grad\-i\-ent of $\,\vt\,$ for
the K\"ah\-ler-sur\-face metric $\,g\,$ on $\,\lz\,$ with
\begin{equation}\label{gvv}
\begin{array}{l}
g(v,v)=g(u,u)=Q,\qquad g(v,u)=g(v,w)=g(u,w)=0,\\
g(w,w')=(\vt\nnh_*\w\nh-\gamma)^{-\nh1}(\vt-\gamma)h(w,w'),
\end{array}
\end{equation}
$w,w'\nh,w''$ always denoting horizontal lifts (that is, projectable
horizontal vector fields), and $\,u\,$ the vector field on $\,\lz\,$ defined 
by $\,u=iv\,$ (multiplication by $\,i\,$ in each fibre). See
\cite[pp.\,1648--1649]{iumj12}. Then $\,g\,$ is a K\"ah\-ler metric for the
al\-most-com\-plex structure $\,J$ obtained by requiring that the
vertical sub\-bundle $\,\vz\hs$ of $\,T\nh M\,$ and $\,\hz\,$ be
$\,J$-in\-var\-i\-ant and the restriction of $\,J\,$ to $\,\vz\nh$, or to
$\,\hz$, coincide with the complex structure of the fibres or, respectively,
with the pull\-back of the complex structure of $\,\varSigma$.

Let $\,M\,$ be the $\,\bbCP^1$ bundle over $\,\varSigma\,$ resulting from the
projective compactification of $\,\lz$. 
According to \cite[Theorem 5.3]{iumj12} and the text preceding it in
\cite{iumj12}, the above construction gives rise to a compact K\"ah\-ler 
surface $\,(M\nh,g)\,$ with the nonconstant trans\-nor\-mal function $\,\vt\,$
having the hol\-o\-mor\-phic gradient $\,v$, and $\,\vt\,$ is 
iso\-par\-a\-met\-ric if and only if $\,\gamma\,$ is constant, while,
conversely, any compact K\"ah\-ler surface carrying a
non-iso\-par\-a\-met\-ric trans\-nor\-mal function with a hol\-o\-mor\-phic
gradient arises from this construction, for a nonconstant $\,\gamma$.

The second paragraph of \cite[Remark 5.2]{iumj12} points out that we can relax
conditions (iii) and (iv), while keeping (ii) and (v) -- (vii), 
so that $\,\varSigma\,$ need not be compact, and $\,Q\,$ is defined and
positive on an open interval. The construction then yields
\begin{equation}\label{qdr}
\mathrm{a\ quadruple\ }\,M\nh,g,J,\vt\,\mathrm{\ with\ the\ same\
properties}
\end{equation}
except compactness of $\,M\nh$,
where $\,M\,$ now is any connected component of the open set in
$\,\lz\smallsetminus\varSigma\,$ defined by requiring that
$\,\vt\ne\gamma\ne\vt\nnh_*\w$ and that the values 
of the norm function $\,r\hs$ lie in the resulting new range.

The classification result of \cite[Theorem 5.3]{iumj12} remains valid, {\it
mutatis mutandis}, in this more general setting: any point at which
$\,d\vt\wedge d\Delta\vt\ne0$, for a trans\-nor\-mal function $\,\vt\,$
with a hol\-o\-mor\-phic gradient on a K\"ah\-ler surface, has a neighborhood
bi\-hol\-o\-mor\-phic\-al\-ly isometric to a noncompact K\"ah\-ler surface
obtained as described in the last paragraph, with nonconstant $\,\gamma$.
To see this, note that, instead of creating the data (i) -- (vii) ``globally''
as in \cite[Sect.\,11]{iumj12}, we may invoke the arguments of
\cite[Sect.\,7]{plms03}, since they all remain valid if $\,\vt\,$ is just
assumed trans\-nor\-mal, rather than iso\-par\-a\-met\-ric.

We now proceed to show that the case of nonconstant $\,\gamma\,$ just mentioned
cannot occur when the resulting K\"ah\-ler surface is weakly Ein\-stein, by
assuming the weakly Ein\-stein property with nonconstant $\,\gamma\,$ and
deriving a contradiction.

This will clearly prove Theorem~\ref{triso}.

To simplify some expressions later, let us note that
\begin{equation}\label{spl}
\begin{array}{l}
h(w,w')\hh D\gamma-h(D\gamma,w')w=-h(J\nnh\nh D\gamma,w)J\hn w'\nh,\\
h(D\gamma,w')J\hn w-h(J\hn w,w')D\gamma=h(D\gamma,w)J\hn w'\nh,\\
h(D\gamma,w')w-h(D\gamma,w)w'\nh=-h(J\hn w,w')J\nnh\nh D\gamma,\\
h(D\gamma,w')J\hn w-h(D\gamma,w)J\hn w'\nh=h(J\hn w,w')D\gamma.
\end{array}
\end{equation}
In fact, the second equality in (\ref{spl}) holds trivially when $\,w,w'$ are
linearly
dependent, and its two sides yield the same $\,h$-in\-ner product with $\,w$,
as well as with $\,w'\nh$. The first (or, third, or fourth) equality arises
from the
second one by replacing $\,w\,$ with $\,J\hn w\,$ (or, applying $\,-J$ to both
sides or, respectively, moving two terms to the other side). 

The third line of (\ref{spl}) gives 
$\,h(D\gamma,D\gamma)w-h(D\gamma,w)D\gamma=h(J\nnh\nh D\gamma,w)J\nnh\nh D\gamma\,$ if
one sets $\,w'\nh=D\gamma$. Applying $\,h(\,\cdot\,,w)\,$ we get
\begin{equation}\label{sth}
h(D\gamma,D\gamma)h(w,w)\,-\,[h(D\gamma,w)]^2\hs=\,[h(J\nnh\nh D\gamma,w)]^2\nh.
\end{equation}
This will explain a seemingly strange presence of
$\,[h(J\nnh\nh D\gamma,w)]^2$ in a later discussion, 
instead of the (expected) left-hand side of (\ref{sth}).

For the 
Le\-vi-Ci\-vi\-ta connection of $\,g\,$ by $\,\nabla\nh$,
\begin{equation}\label{lcc}
\begin{array}{l}
\nabla\nh_{\!v}\w v=-\nabla\nh_{\!u}\w u=\psi\hh v,\qquad
\nabla\nh_{\!v}\w u=\nabla\nh_{\!u}\w v=\psi\hh u,\\
\nabla\nh_{\!v}\w w=\nabla\nh_{\!w}\w v=\phi\hh w,\qquad
\nabla\nh_{\!u}\w w=\nabla\nh_{\!w}\w u=\phi\hh J\hn w,\\
\nabla\nh_{\!w}\w w'\nh=D\hn_{\!w}\w w'\nh
+S\hh[h(D\gamma,w)w'\nh+h(J\nnh\nh D\gamma,w)J\hn w']\\
\hskip30pt-\,[2(\vt\nnh_*\w\nh-\gamma)]^{-\nh1}\hs[h(w,w')\hh v
+h(J\hn w,w')\hh u].
\end{array}
\end{equation}
Here $\,D\,$ stands both for the Le\-vi-Ci\-vi\-ta connection of the
base-sur\-face metric $\,h\,$ and for the $\,h$-grad\-i\-ent,
while $\,\psi,\phi,S\,$ are the functions given by
\begin{equation}\label{tps}
2\psi=\dot Q,\qquad2\phi=Q/(\vt-\gamma),\qquad
2S=(\vt\nnh_*\w\nh-\gamma)^{-\nh1}\nh-(\vt-\gamma)^{-\nh1}\nh,
\end{equation}
where $\,(\,\,)\dot{\,}=\,d/d\vt$. We have used the first line of
(\ref{spl}) to replace the expression 
$\,S\hh[h(D\gamma,w)w'\nh+h(D\gamma,w')w-h(w,w')D\gamma]\,$ appearing (with
a slightly different notation) in the description of
$\,\nabla$ \cite[p.\,1648]{iumj12}, by
$\,S\hh[h(D\gamma,w)w'\nh+h(J\nnh\nh D\gamma,w)J\hn w']$.

Applying the equality $\,2\hh\ric(v,\,\cdot\,)=-d\Delta\vt\,$ in
Remark~\ref{calab} to $\,\vt\,$ and $\,v$ appearing above, that is,
in \cite[formula\,(5.2)]{iumj12}, and noting that $\,\Delta\vt=2(\psi+\phi)$, 
while $\,d_v\w$ acts on functions of $\,\vt\,$ as $\,Q\,d/d\vt$, we get
\begin{equation}\label{ric}
\begin{array}{l}
2(\vt-\gamma)^2\hh\ric(v,v)=2(\vt-\gamma)^2\hh\ric(u,u)
=QP\nh,\\
\mathrm{for\ }\,P\nh=Q-(\vt-\gamma)\dot Q-(\vt-\gamma)^2\nh\ddot Q,\\
\ric(v,u)=0,\quad2(\vt-\gamma)^2\hh\ric(v,w)=-Q\hh h(D\gamma,w),\\
2(\vt-\gamma)^2\hh\ric(u,w)=-Q\hh h(J\nnh\nh D\gamma,w),
\end{array}
\end{equation}
using $\,J$-in\-var\-i\-ance of $\,\ric\,$ to derive the equations involving
$\,u\,$ from those for $\,v$. This shortcut describes all components of
$\,\ric\,$ except $\,\ric(w,w')\,$ for two horizontal vectors $\,w,w'\nh$.

With the sign convention
$\,R\hs(v,w)u=\nabla\!_{[v,w]}\w u+\nabla\!_{w}\w\nabla\!_{v}\w u-
\nabla\!_{v}\w\nabla\!_{w}\w u$, (\ref{lcc}) and the equality
$\,-Q\dot\phi=2(\phi-\psi)\phi$, immediate from (\ref{tps}),
yield the following equalities, describing all components of the
$\,(0,4)\,$ curvature tensor except those involving four
horizontal vectors:
\[
\begin{array}{l}
2R(v,u)v=-Q\ddot Q\hh u,\quad 2R(v,u)u
=Q\ddot Q\hh v,\quad R(v,u)w=2(\phi-\psi)\phi\hh J\hn w,\\
R(v,w)v=(\phi-\psi)\phi\hh w,\quad
R(v,w)u=(\phi-\psi)\phi\hh J\hn w,\\
2R(v,w)w'\nh=(\vt\nnh_*\w\nh-\gamma)^{-\nh1}(\psi-\phi)[h(w,w')\hh v+h(J\hn w,w')\hh u]
-(\vt-\gamma)^{-\nh2}Q\hh h(J\nnh\nh D\gamma,w)J\hn w'\nh,\\
R(u,w)v=(\psi-\phi)\phi\hh J\hn w,\quad
R(u,w)u=(\phi-\psi)\phi\hh w,\\
2R(u,w)w'\nh=(\vt\nnh_*\w\nh-\gamma)^{-\nh1}(\phi-\psi)[h(J\hn w,w')\hh v
-h(w,w')\hh u]+(\vt-\gamma)^{-\nh2}Q\hh h(D\gamma,w)J\hn w'\nh,\\
2R(w,w')v=-(\vt-\gamma)^{-\nh2}Q\hh h(J\hn w,w')J\nnh\nh D\gamma
+2(\vt\nnh_*\nh-\gamma)^{-\nh1}(\phi-\psi)h(J\hn w,w')u,\\
2R(w,w')u\nh=\nnh(\vt-\gamma)^{-\nh2}Q\hh h(J\hn w,w')D\gamma
+\nh2(\vt\nnh_*\nh-\gamma)^{-\nh1}(\psi-\phi)h(J\hn w,w')v.
\end{array}
\]
where we simplified $\,2R(w,w')v\,$ using the third line of (\ref{spl}) and
then derived the expression for $\,2R(w,w')u\,$ from the fact that
$\,R(w,w')u=J[R(w,w')v]$. 

For dimensional reasons, a ``horizontal'' component
$\,R(w,w'\nh,w''\nh,w''')\,$ must be given by
\begin{equation}\label{mbe}
2(\vt\nnh_*\nh-\gamma)(\vt-\gamma)\hh R(w,w'\nh,w''\nh,w''')\,
=\,Z\hh[h(w,w'')g(w'\nh,w''')\,-\,h(w',w'')g(w,w''')]
\end{equation}
for some function $\,Z\,$ not depending
on $\,w,w'\nh,w''\nh,w'''\nh$. Before determining what $\,Z\,$ is, we  
now characterize the weak\-ly-Ein\-stein case as a condition imposed
on $\,Z$. First,
\begin{equation}\label{rvp}
\begin{array}{l}
2(\vt\nnh_*\nh-\gamma)(\vt-\gamma)\hh R(w,w')w''\hs\,
=\,Z\hh[h(w,w'')\hh w'\hs-\,h(w',w'')\hh w]\\
\hskip130pt+\,\,h(J\hn w,w')h(J\nnh\nh D\gamma,w'')\hh v
-h(J\hn w,w')h(D\gamma,w'')\hh u,
\end{array}
\end{equation}
which now easily implies that, with 
$\,\widetilde Z=Z+2(\vt-\gamma)(\phi-\psi)$,
\[
2(\vt\nnh_*\nh-\gamma)(\vt-\gamma)\hs\ric(w,w')\,
=\,\widetilde Zh(w,w').
\]
The Ric\-ci en\-do\-mor\-phism of $\,T\nh M\,$ acts by
\begin{equation}\label{ren}
\begin{array}{l}
2(\vt-\gamma)^3\ric\hh v=(\vt-\gamma)\hh Pv-(\vt\nnh_*\nh-\gamma)\hh QD\gamma,
\quad \ric\hh u=J\ric\hh v,\\
2(\vt-\gamma)^2\ric\hh w=\widetilde Z\hh w-h(D\gamma,w)v
-h(J\nnh\nh D\gamma,w)u.
\end{array}
\end{equation}
The value assigned by $\,R\hh\ric\,$ to each of the six pairs
\[
(v,v),\,\,\,(v,u),\,\,\,(v,w),\,\,\,(u,u),\,\,\,(u,w),\,\,\,(w,w')
\]
of vector fields equals the trace of the composition in which 
the Ric\-ci en\-do\-mor\-phism, with (\ref{ren}), is followed by
\begin{equation}\label{flw}
R(v,\,\cdot\,)v,\,\,\,R(v,\,\cdot\,)u,\,\,\,R(v,\,\cdot\,)w,\,\,\,
R(u,\,\cdot\,)u,\,\,\,R(u,\,\cdot\,)w,\,\,\,R(w,\,\cdot\,)w'.
\end{equation}
Due to Her\-mit\-i\-an symmetry of $\,R\hh\ric\,$ (that is, its
$\,J$-in\-var\-i\-ance) we only need to consider three of the six pairs:
$\,(v,v),\,(v,w)\,$ and $\,(w,w')$.
First,
\begin{equation}\label{rvs}
\begin{array}{l}
R(v,\,\cdot\,)v\,\mathrm{\ \ and\ \ }\,R(u,\,\cdot\,)u\,\mathrm{\ \ send\ the\
triple\ \ }\,(v,u,w)\,\mathrm{\ \ to}\\ 
(0,-Q\ddot Q\hh u/2,(\phi-\psi)\phi w)\,\mathrm{\ and\
}\,(-Q\ddot Q\hh v/2,0,(\phi-\psi)\phi w).
\end{array}
\end{equation}
Hence, with
$\,\widetilde Z=Z+2(\vt-\gamma)(\phi-\psi)\,$ as before,
\[
4(\vt-\gamma)^2[R\hh\ric](v,v)
=4(\phi-\psi)\phi\hh\widetilde Z-Q\ddot QP\nh.
\]
Similarly,
\[
4(\vt-\gamma)^4[R\hh\ric](v,w)
=[4(\vt-\gamma)^2(\psi-\phi)\phi
-Q\widetilde Z]\hh h(D\gamma,w),
\]
and $\,4(\vt\nnh_*\nh-\gamma)(\vt-\gamma)^4[R\hh\ric](w,w')\,$ equals
$\,h(w,w')\,$ times
\begin{equation}\label{hwt}
(\vt-\gamma)Z\hn\widetilde Z
+2(\vt\nnh_*\nh-\gamma)\hh Qh(D\gamma,\nnh D\gamma)
+2(\vt-\gamma)^2(\phi-\psi)P
\end{equation}
or, equivalently, $\,4(\vt-\gamma)^5[R\hh\ric](w,w')\,$ equals 
$\,g(w,w')\,$ times (\ref{hwt}). In the case where $\,\gamma\,$ is
nonconstant, for $\,R\hh\ric\,$ to be a functional multiple of $\,g\,$ it is
necessary and sufficient that 
$\,4(\vt-\gamma)^2(\psi-\phi)\phi=Q\widetilde Z\,$ and that 
$\,(\vt-\gamma)^3\nh[4(\phi-\psi)\phi\hh\widetilde Z-Q\ddot QP]\,$ be equal to
$\,Q\,$ times (\ref{hwt}).
Since $\,2(\vt-\gamma)\hh\phi=Q$, the first condition amounts to
\begin{equation}\label{fcn}
\widetilde Z=2(\vt-\gamma)(\psi-\phi),\quad\mathrm{that\ is,\
}\,Z=4(\vt-\gamma)(\psi-\phi).
\end{equation}
Assuming (\ref{fcn}), and noting that
\[
-(\vt-\gamma)^2\ddot Q=P-Q+(\vt-\gamma)\dot Q=P\nh+2(\vt-\gamma)(\psi-\phi),
\]
we rewrite the second condition as
\begin{equation}\label{scn}
2(\vt\nnh_*\nh-\gamma)\hh Qh(D\gamma,\nnh D\gamma)
=(\vt-\gamma)[P\nh+6(\vt-\gamma)(\psi-\phi)][P\nh-2(\vt-\gamma)(\psi-\phi)].
\end{equation}
Since $\,P\nh=Q-(\vt-\gamma)\dot Q-(\vt-\gamma)^2\nh\ddot Q\,$ while, by
(\ref{tps}), $\,2\psi=\dot Q\,$ and $\,2(\vt-\gamma)\phi=Q$, the two factors
in square brackets on the right-hand side of (\ref{scn}) are equal to
\[
-(\vt-\gamma)^2\nh\ddot Q\,+\,2[(\vt-\gamma)\dot Q-Q]\,\mathrm{\ \ and\ \
}\,-\nnh(\vt-\gamma)^2\nh\ddot Q\,-\,2[(\vt-\gamma)\dot Q-Q],
\]
so that their product is
$\,(\vt-\gamma)^4\nh\ddot Q^2\nh-4[(\vt-\gamma)\dot Q-Q]^2$ and we may rewrite
the right-hand side of (\ref{scn}), divided by $\,Q$, as
\[
\begin{array}{l}
-Q^{-\nh1}\nh\ddot Q^2\nh\gamma^5\nh+5\vt Q^{-\nh1}\nh\ddot Q^2\nh\gamma^4\nh
+2Q^{-\nh1}\nh[2\dot Q^2\nh-5\vt^2\nh\ddot Q^2]\gamma^3\nh
+2[4\dot Q+6\vt Q^{-\nh1}\nh\dot Q^2\nh
+5\vt^3\nh Q^{-\nh1}\nh\ddot Q^2]\gamma^2\\
\hskip20pt+[4Q-16\vt\dot Q+12\vt^2\nh Q^{-\nh1}\nh\dot Q^2\nh
-5\vt^4\hn Q^{-\nh1}\nh\ddot Q^2]\gamma
-4\vt Q+8\vt^2\nh\dot Q-4\vt^3\nh Q^{-\nh1}\nh\dot Q^2\nh
+\vt^5\nh Q^{-\nh1}\nh\ddot Q^2\nh.
\end{array}
\]
This is a quintic polynomial in $\,\gamma$, with coefficients that are
functions of $\,\vt$, which equals -- according to (\ref{scn}) -- a function
on base-sur\-face $\,\varSigma$. Applying $\,d_v\w$ to the latter function
we get $\,0\,$ and, since $\,d_v\w$ acts on functions of $\,\vt\,$
as $\,Q\,d/d\vt$, we conclude that {\it the coefficients of the above 
quintic polynomial are constant functions of\/} $\,\vt$. However, constancy of
both $\,Q^{-\nh1}\nh\ddot Q\,$ and $\,\vt Q^{-\nh1}\nh\ddot Q\,$ means that
$\,\ddot Q=0$, and hence $\,\dot Q\,$ is constant. Looking at the
coefficients of $\,\gamma^3$ and $\,\gamma\,$ we now see that $\,Q\,$ must be
constant, and hence zero.

This contradiction proves that $\,\gamma\,$ is constant. In other words, our
$\,\vt$, 
besides being trans\-nor\-mal, is also iso\-par\-a\-met\-ric.

\section{The lo\-cal-struc\-ture theorem}\label{ls}
By a {\it special K\"ah\-ler-Ric\-ci potential\/} 
\cite[Sect.\,7]{plms03} 
on a K\"ah\-ler manifold $\,(M\nh,g,J)\,$ one means any nonconstant function
$\,\vt\,$ on $\,M$ having a 
real-hol\-o\-mor\-phic gradient 
$\,v=\nabla\nh\vt$ 
for which, at points where $\,v\ne0$, all 
nonzero vectors orthogonal to $\,v\,$ and $\,J\nh v\,$ are eigen\-vec\-tors 
of both $\,\nabla\nh d\vt$ and the Ric\-ci tensor $\,\ric$. Such 
quadruples $\,(M\nh,g,J,\vt)\,$ have been completely described, both locally 
\cite{plms03} and in the compact case \cite{jram06}.

In the case of K\"ah\-ler surfaces, $\,\vt\,$ as above is nothing else than
a nonconstant iso\-par\-a\-met\-ric function with a hol\-o\-mor\-phic
gradient. In fact, generally, on any Riemannian manifold,
for $\,v=\nabla\nh\vt\,$ and $\,Q=g(v,v)\,$ one has
$\,d\hh Q=2[\nabla\nh d\vt](v,\,\cdot\,)\,$ while, if $\,v\,$ is
real hol\-o\-mor\-phic on a K\"ah\-ler manifold, 
$\,2\hh\ric(v,\,\cdot\,)=-d\Delta\vt\,$ (see Remark~\ref{calab}).

The construction in \cite{plms03} is a special case of
\cite[pp.\,1648--1649]{iumj12}, with constant $\,\gamma$. The nonzero constant 
$\,\vt\nnh_*\w\nh-\gamma\,$ in (\ref{gvv}) is replaced with $\,\ve/2$, where
$\,\ve=\pm1\,$ (with no loss of generality, since the base-sur\-face metric
$\,h\,$ can be rescaled). Now, from (\ref{gvv}), as in
\cite[pp.\,791--792 and Sect.\,16]{plms03}
\begin{equation}\label{gvc}
\begin{array}{l}
g(v,v)=g(u,u)=Q,\qquad g(v,u)=g(v,w)=g(u,w)=0,\\
g(w,w')=2\ve(\vt-\gamma)h(w,w'),
\end{array}
\end{equation}
$w,w'$ still denoting horizontal lifts (projectable
horizontal vector fields). In terms of the vertical and horizontal
distributions $\,\mathcal{V}=\mathrm{Span}\hs(v,u)\,$ and 
$\,\mathcal{H}=\mathcal{V}^\perp\nnh$,
\begin{equation}\label{rce}
\begin{array}{l}
\ric=\mu g\,\mathrm{\ on\ }\,\mathcal{V}\nnh,\nnh\qquad\ric
=\lambda g\,\mathrm{\ on\
}\,\mathcal{H},\nnh\qquad\ric(\mathcal{V}\nnh,\mathcal{H})=\{0\},\nnh\quad
\mathrm{where}\\
\mu=-\dot Y/2,\quad\lambda
=[2\ve(\vt-\gamma)]^{-\nh1}(K-\ve Y),\mathrm{\ for\ }\,Y\nh=2(\psi+\phi)
\end{array}
\end{equation}
(so that $\,Y\nnh=\Delta\vt$) and $\,K\,$ is the Gauss\-i\-an curvature of 
$\,h$. See (\ref{ric}) and 
\cite[formula\,(7.4), the lines following (8.1), and (b) in
Sect.\,16]{plms03}, where our $\,\gamma\,$ is denoted by $\,c$. Thus,
once we identify $\,v\,$ with $\,g(v,\,\cdot\,)$, and similarly for 
$\,u$,
\begin{equation}\label{rre}
\begin{array}{l}
\hs\ric\,\,=\,\,\lambda g\,+\,(\mu-\lambda)Q^{-\nh1}\nh(v\otimes v\,
+\,u\otimes u)\mathrm{,\ \ and\ \,hence}\\  
R\hh\ric\,=\,\lambda\hs\ric\,
+\,(\mu-\lambda)Q^{-\nh1}\nh[R(v,\,\cdot\,,v,\,\cdot\,)
+R(u,\,\cdot\,,u,\,\cdot\,)].
\end{array}
\end{equation}
From (\ref{rvs}) and (\ref{rce}), $\,R\hh\ric\,$ treated as an 
en\-do\-mor\-phism of $\,T\nh M\,$ sends $\,v,u,w\,$ to
\[
[\lambda\mu+(\lambda-\mu)\ddot Q/2]v,\quad
[\lambda\mu+(\lambda-\mu)\ddot Q/2]u,\quad
[\lambda^2\nh+(\vt-\gamma)^{-\nh1}\nnh(\mu-\lambda)(\phi-\psi)]w,
\]
$\,\psi,\phi\,$ being -- as in (\ref{tps}) -- the functions given by
\begin{equation}\label{psf}
2\psi=\dot Q,\qquad\quad2\phi=Q/(\vt-\gamma)
\end{equation}
Thus,
$\,R\hh\ric\,$ is a functional multiple of the identity (cf.\
Theorem~\ref{equiv}) if and only if
\begin{equation}\label{eul}
(\vt-\gamma)^2\nh\ddot Q\,+\,2Q\,=\,\ve K(\vt-\gamma),
\end{equation}
as long as we exclude the Ein\-stein case by assuming that $\,\mu\ne\lambda$.
(See Remark~\ref{divby} below.)

One immediate conclusion is that $\,K$, the Gauss\-i\-an curvature of the
base-sur\-face metric $\,h$, must be constant: $\,d_w\w\vt=0\,$ and $\,Q\,$ is
a function of $\,\vt$, so $\,d_w\w$ applied to (\ref{eul}) yields
$\,d_w\w K=0$. Solving (\ref{eul}), we see that
\begin{equation}\label{qeq}
\begin{array}{l}
Q\,\mathrm{\hn\ as\hn\ a\hn\ function\hn\ of\hn\ }\,\vt\,\mathrm{\hn\ equals\
}\,\ve K(\vt-\gamma)/2\,\mathrm{\hn\ plus\hn\ a\hn\ linear\hn\
combination\hn\ of}\\
|\vt-\gamma|^{1/2}\cos\hs[\nh\sqrt{7\,}(\log|\vt-\gamma|)/2]\hs\,\,\mathrm{\
and\ }\,\,\hs|\vt-\gamma|^{1/2}\sin\hs[\nh\sqrt{7\,}(\log|\vt-\gamma|)/2].
\end{array}
\end{equation}
\begin{thm}\label{weskr}Defining the quadruple\/ $\,(M\nh,g,J,\vt)\,$ by the
local version\/ {\rm(\ref{qdr})} of the construction in
Section\/~{\rm\ref{pm}}, for\/ $\,Q>0\,$ as in\/ {\rm(\ref{qeq})}, we obtain 
a nonconstant iso\-par\-a\-met\-ric function\/ $\,\vt\,$ with a
real-hol\-o\-mor\-phic gradient on a weakly Ein\-stein K\"ah\-ler surface.

Conversely, up to local bi\-hol\-o\-mor\-phic isometries, any
nonconstant trans\-nor\-mal function with a 
real-hol\-o\-mor\-phic gradient on a weakly Ein\-stein K\"ah\-ler surface
arises in the manner just described.
\end{thm}
\begin{pf}The first part of the theorem is immediate from the
preceding discussion and
Theorem~\ref{equiv}(b). 
The final clause is, according to the second paragraph of this section, 
a special case of \cite[Theorem 18.1]{plms03}, since 
`trans\-nor\-mal' in our situation implies `iso\-par\-a\-met\-ric'
as a consequence of Theorem~\ref{triso}.\qed
\end{pf}
\begin{rem}\label{divby}We arrived at (\ref{eul}) after dividing an
intermediate equality by $\,\lambda-\mu\,$ (``excluding the Ein\-stein 
case''). 
By (\ref{rce}) and (\ref{psf}),
$\,2(\vt-\gamma)(\lambda-\mu)
=(\vt-\gamma)^2\nh\ddot Q-2Q+\ve K(\vt-\gamma)$,
so that $\,g\,$ is an Ein\-stein metric on any nonempty open set on which
\begin{equation}\label{lmm}
(\vt-\gamma)^2\nh\ddot Q\,-\,2Q\,=\,-\ve K(\vt-\gamma).
\end{equation}
Dividing by $\,\lambda-\mu\,$ in the non-Ein\-stein case is allowed, without
assuming real-an\-a\-lyt\-ic\-i\-ty, even 
when such a nonempty open set exists. In fact, for some open interval 
$\,I\hh'$ of the variable $\,\vt$, one then has (\ref{lmm}) on a proper subset
$\,P\hs$ of $\,I\hh'$ with a nonempty interior (which is a disjoint union of
a countable family of open intervals), and on $\,I\hh'\nh\smallsetminus P\nh$,
which is also such a union, (\ref{eul}) holds. If
$\,I\hh''\nh\subseteq P\nh$, or
$\,I\hh''\nh\subseteq I\hh'\nh\smallsetminus P\nh$, 
is one of those countably many subintervals, then an endpoint 
$\,\vt_0\w\in I\hh'$ of $\,I\hh''$ is a cluster point for both 
countable unions. The Taylor series of $\,Q\,$ at 
$\,\vt_0\w$ thus satisfies the series versions of both equations (\ref{eul})
and (\ref{lmm}), which determines it uniquely -- see below -- while $\,Q\,$
is a solution to (\ref{lmm}) or, respectively, (\ref{eul}), on the half-open
interval $\,I\hh''\nnh\cup\{\vt_0\w\}$. Namely, treated as a formal power
series, our Taylor series, satisfying both (\ref{eul}) and (\ref{lmm}), 
must correspond to
\begin{equation}\label{qee}
Q\,=\,\ve K(\vt-\gamma)/2\hh,
\end{equation}
and consequently have $\,Q(\vt_0\w)=\ve K(\vt_0\w-\gamma)/2\,$ and 
$\,\dot Q(\vt_0\w)=\ve K/2$. Our function $\,Q\,$ is thus given by
(\ref{qee}), due to uniqueness of a solution to (\ref{lmm}), or
(\ref{eul}), on $\,I\hh''\nnh\cup\{\vt_0\w\}\,$ with the above initial data. 
We therefore have (\ref{qee}) on an open dense set, and hence everywhere, 
in $\,I\hh'\nnh$. In other words, (\ref{eul}) {\it always follows in the 
weak\-ly-Ein\-stein non-Ein\-stein case}, since otherwise, as we just saw, 
$\,Q=\ve K(\vt-\gamma)/2$, implying the Ein\-stein property -- namely, 
the weak\-ly-Ein\-stein metric resulting from the 
construction based on (\ref{eul}), and mentioned in Theorem~\ref{weskr}, is 
Ein\-stein if and only if $\,Q=\ve K(\vt-\gamma)/2$.
\end{rem}
\begin{rem}\label{notlh}The weak\-ly-Ein\-stein metric arising in 
Theorem~\ref{weskr} are not locally homogenous case except in the Ein\-stein
case (when $\,Q=\ve K(\vt-\gamma)/2$, cf.\ Remark~\ref{divby}). In fact,
for the eigen\-value function $\,\mu\,$ of the Ric\-ci tensor, (\ref{rce})
gives $\,-\nh2\mu(\vt-\gamma)^2\nh=(\vt-\gamma)^2\nh\ddot Q+
(\vt-\gamma)\dot Q-Q\,$ which, by (\ref{eul}), equals
$\,(\vt-\gamma)\dot Q+\ve K(\vt-\gamma)-3Q$. Solving the resulting equation
$\,(\vt-\gamma)\dot Q-3Q=-\ve K(\vt-\gamma)-\nh2\mu(\vt-\gamma)^2$ under the
assumption that$\,\mu\,$ is constant, we get
$\,Q=a(\vt-\gamma)^3\nh+\mu(\vt-\gamma)^2\nh+\ve K(\vt-\gamma)/2$ with a
constant $\,a$, which is not of the form (\ref{qeq}) except when
$\,a=\mu=0\,$ and $\,Q=\ve K(\vt-\gamma)/2$.
\end{rem}
\begin{rem}\label{prdct}In \cite[sixth line of Sect.\,8]{plms03} there is
a second case, with the factor $\,2\ve(\vt-\gamma)$ in (\ref{gvc}) 
replaced by $\,1$. We leave it out of our discussion, since this is precisely
the case of product metrics \cite[Corollary 13.2 and (c) in Sect.\,16]{plms03}.
See also Remark~\ref{justf}.
\end{rem}

\section{Nonrealizable boundary conditions}\label{br}
The following lemma is a crucial step in the proof of Theorem~\ref{trhol}.
\begin{lem}\label{nexst}For\/ $\,F:\bbR\to\bbR\,$ given by\/
$\,F(\alpha)=e^{-\alpha\cot c}\sin\alpha$, where\/
$\,c=\tan^{-\nh1}\hskip-3pt\sqrt{7\,}$, there do not exist\/
$\,\alpha,\beta\in\bbR\,$ such that\/ $\,\alpha<\nh\beta\,$ and\/
$\,F(\alpha)-F(\beta)=F\hh'\nnh(\alpha)+F\hh'\nnh(\beta)=0$, while\/
$\,F\,$ is nonzero everywhere in the open interval\/ $\,(\alpha,\beta)$.
\end{lem}
\begin{pf}Since $\,F\,$ is pos\-i\-tive/neg\-a\-tive on the interval
$\,(n\pi,(n+1)\pi)\,$ for each even/odd integer $\,n$, such $\,\alpha,\beta$,
if they existed, would both lie in one of the closed intervals 
$\,I\nnh_n\w=[n\pi,(n+1)\pi]$. The identity 
$\,F(\alpha+\pi)=-e^{-\pi\cot c}\nnh F(\alpha)\,$ shows that the graph of
$\,F\,$ is the same, up to vertical rescaling, on all $\,I\nnh_n\w$, and so
it suffices to consider the case $\,n=1$, that is, prove the nonexistence of
$\,\alpha,\beta\in[\hh0,\pi]\,$ with the stated properties.

Our $\,F\,$ assumes in $\,[\hh0,\pi]\,$ the minimum value 
$\,0$, just at the endpoints, and the maximum value $\,F(c)$, only 
at $\,c$, as induction on $\,q\ge0\,$ gives, for $\,F^{(q)}\nh
=d\hh^q\nnh\hn F\nnh/\nh d\alpha^q\nh$,
\begin{equation}\label{fpr}
\begin{array}{rl}
\mathrm{i)}&F^{(q)}\nnh(\alpha)\sin^q\nnh c\hh\,
=\hh\,(-\nh1)^q\nh e^{-\alpha\cot c}\sin\hs(\alpha-qc),\\
\mathrm{ii)}&F'\hn\nnh>0\,\,\mathrm{\ on\ }\,[\hh0,c)\mathrm{,\ while\ 
}\hs\,F'\hn\nnh<0\,\hs\mathrm{\ on\ }\,(c,\pi].
\end{array}
\end{equation}
Numerically, $\,c\approx1.209429$, and so $\,c\in(0,\pi/2)$. 
By (\ref{fpr}-ii), $\,F\,$ maps both $\,[\hh0,c)\,$ and 
$\,(c,\pi]\,$ dif\-feo\-mor\-phic\-al\-ly onto $\,[\hh0,F(c))$. Thus, for
every $\,\alpha\in[\hh0,\pi]\smallsetminus\{c\}\,$  there exists a unique 
$\,\beta=\beta(\alpha)\in[\hh0,\pi]\smallsetminus\{c\}\,$ with
$\,\beta\ne\alpha\,$ and $\,F(\beta)=F(\alpha)$. Now
\begin{equation}\label{dsd}
\begin{array}{l}
\mathrm{setting\ }\,\beta(c)=c\,\mathrm{\ we\ obtain\ a\ decreasing\
dif\-feo\-mor}\hyp\\
\mathrm{phism\
}\,[\hh0,\pi]\ni\alpha\mapsto\beta(\alpha)\in[\hh0,\pi]\,\mathrm{\ with\
}\,\beta'\nnh(c)=-\nh1.
\end{array}
\end{equation}
Namely, smoothness at $\,c\,$ in (\ref{dsd}) follows from the Morse 
lemma: $\,(c,c)\,$ is a nondegenerate critical point with the index $\,0\,$ 
and value $\,0\,$ for the function on 
$\,\bbR\nh^2$ sending $\,(\alpha,\beta)$ to $\,F(\alpha)-F(\beta)$, and so
the zeros of this function near $\,(c,c)\,$ form two smooth curves intersecting
transvesally at $\,(c,c)$. The two curves are the graphs of the identity
function and the function $\,\alpha\mapsto\beta(\alpha)\,$ in (\ref{dsd}). 
Hence $\,\beta'\nnh(c)\,$ exists and equals 
$\,-\nh1$, as the invariance 
of the graph in (\ref{dsd}) under the switch of $\,\alpha\,$ with $\,\beta\,$
causes the vector $\,(1,-\nh1)\,$ to be tangent to the graph at $\,(c,c)$,
and thus ensuring that the graph does not pass through $\,(c,c)\,$ vertically. 
Now (\ref{fpr}-i) gives
$\,F\hh''\nnh(\alpha)\sin^2c=e^{-\alpha\cot c}\sin\hs(\alpha-2c)\,$ and
$\,F\hh'''\nnh(\alpha)\sin^3c=e^{-\alpha\cot c}\sin\hs(3c-\alpha)$,
so that,
\begin{equation}\label{cfd}
\begin{array}{l}
\mathrm{as\ }\,\,3c\,-\,\pi\,\approx\,0.4867\,\,\mathrm{\ lies\ in\
}\,\,(0,\pi/4)\mathrm{,\ and\ hence\ }\,\,0\hs<3c-\pi\,<\hs c,\\
F\hh''\nh\mathrm{\ decreases\ on\ }\,[\hh0,3c-\pi]\,\mathrm{\ \,from\ 
}-\hskip-3.5pt2\cot c\,\mathrm{\ \ to\ }\,F\hh''\nnh(3c-\pi)\mathrm{,\
increas}\hyp\\
\mathrm{es\ on\ }\,\hs[3c-\pi,2c]\,\mathrm{\hs\ from\
}\hs\,F\hh''\nnh(3c-\pi)\,\,\mathrm{\ to\ }\hs\,0\mathrm{,\ and\ increases\ on\
}\,\hs[2c,\pi].
\end{array}
\end{equation}
The identity $\,F(\alpha)=F(\beta(\alpha))$ and the chain rule give
\begin{equation}\label{bpa}
\beta'\nnh(\alpha)\,
=\,\frac{F\hh'\nnh(\alpha)}{F\hh'\nnh(\beta(\alpha))}\,\mathrm{\ whenever\ 
}\,\alpha\in[\hh0,c),
\end{equation}
and so the assertion of the lemma amounts to
\begin{equation}\label{bpb}
\beta'\nnh(\alpha)\,\ne\,-\nh1\,\mathrm{\ for\ all\ }\,\alpha\in[\hh0,c).
\end{equation}
We will now obtain (\ref{bpb}), and complete the proof, by showing that
\begin{equation}\label{bpm}
\beta'\,\,<\,\,-\nh1\,\mathrm{\ on\ the\ interval\ }\,[\hh0,c).
\end{equation}
Let $\,\alpha_0\w=\beta(2c)$. By (\ref{dsd}), $\,\alpha_0\w\in(0,c)$. Also,
\begin{equation}\label{lmo}
\beta'\nnh(0)\,\,<\,\,-\nh1\,\mathrm{\ and\
}\,\beta'\nnh(\alpha_0\w)\,\,<\,\,-\nh1.
\end{equation}
In fact, as
$\,\beta(0)=\pi$, (\ref{bpa}) and (\ref{fpr}-i) with $\,q=1\,$ yield 
$\,\beta'\nnh(0)=-\hh e^{\pi\cot c}\nh<-\nh1\,$ while, 
with the approximate values 
$\,\alpha_0\w\approx0.3017\,$ and
$\,\beta(\alpha_0\w)=2c\approx2.418858\,$ 
provided by {\it Mathematica}, (\ref{bpa}) and (\ref{fpr}-i) for $\,q=1\,$ 
give $\,\beta'\nnh(\alpha_0\w)\approx-\nh1.8755$.

Differentiating the identity 
$\,F(\alpha)=F(\beta(\alpha))\,$ twice we get
\begin{equation}\label{fpb}
F\hh'\nnh(\beta(\alpha))\beta''\nnh(\alpha)
=F\hh''\nnh(\alpha)-F\hh''\nnh(\beta(\alpha))[\beta'\nnh(\alpha)]^2\nh.
\end{equation}
One verifies numerically that 
$\,F\hh''\nnh(0)/F\hh''\nnh(c)=2e^{c\cot c}\cos c\approx1.169>1$, 
and so, by (\ref{cfd}),
\begin{equation}\label{fpp}
F\hh''\nnh(\alpha)<F\hh''\nnh(c)<F\hh''\nnh(\beta)<0\,\mathrm{\ whenever\
}\,0\le\alpha<c<\beta<2c.
\end{equation}
For $\,\alpha\in[0,\alpha_0\w]\,$ we have $\,\beta(\alpha)\ge2c\,$ and, by
(\ref{cfd}), the right-hand side of (\ref{fpb}) is negative, while
$\,F\hh'\nnh(\beta(\alpha))<0\,$ due to (\ref{fpr}-ii). Thus, 
$\,\beta''\nnh(\alpha)>0$, so that $\,\beta'$ is increasing on 
$\,[0,\alpha_0\w]\,$ and, by (\ref{lmo}), $\,\beta'\nh<-\nh1\,$ 
on $\,[0,\alpha_0\w]$.

The remaining part of (\ref{bpm}) is the inequality
$\,\beta'\nh<-\nh1\,$ on $\,(\alpha_0\w,c)$. We establish it by
showing that its negation leads to a contradiction. Namely, suppose that 
$\,\beta'\nnh(\alpha)\ge-\nh1$ for some $\,\alpha\in(\alpha_0\w,c)$.
{\it At any such\/ $\,\alpha$,  negativity of\/ $\,\beta'\nh$, due to\/
{\rm(\ref{dsd})}, gives\/ $\,|\beta'\nnh(\alpha)|\le1$, and so, by\/
{\rm(\ref{fpp})} and\/} {\rm(\ref{fpb})}, $\,\beta''\nnh(\alpha)>0$. (In fact,
$\,\beta(\alpha)\in(c,2c)$, cf.\ (\ref{dsd}), so that (\ref{fpr}-ii) and
(\ref{cfd}) imply negativity of
$\,F\hh'\nnh(\beta(\alpha))$, $\,F\hh''\nnh(\alpha)\,$ and
$\,F\hh''\nnh(\beta(\alpha))$, and (\ref{fpp}) with
$\,1-[\beta'\nnh(\alpha)]^2\nh\ge0\,$ yields
$\,F\hh''\nnh(\alpha)-F\hh''\nnh(\beta(\alpha)[\beta'\nnh(\alpha)]^2\nh
<F\hh''\nnh(\beta(\alpha))-F\hh''\nnh(\beta(\alpha))[\beta'\nnh(\alpha)]^2\nh
\le0$.) Remark~\ref{posdr} applied to $\,\psi=\beta'$
and $\,\lambda=-\nh1$ now shows that $\,\beta'$ has a limit 
at $\,c\,$ greater than $\,-\nh1$. This contradicts the
equality $\,\beta'\nnh(c)=-\nh1\,$ in (\ref{dsd}), completing the proof.\qed
\end{pf}

\section{Proof of Theorem~\ref{trhol}}\label{pl}
Under our assumptions, according to Theorem~\ref{triso}, the second paragraph
of Sect.\,\ref{ls} and \cite[Theorem 16.3]{jram06}, $\,(M\nh,g,J)\,$ arises,
up to bi\-hol\-o\-mor\-phic isometries, from the ``compact'' version
(\ref{qdr}) of the construction in Sect.\,\ref{pm}, using some data (i) --
(vii) with $\,\vt\nnh_*\w\nh-\gamma\,$ in (\ref{gvv}) replaced by
$\,\ve/2$, where $\,\ve=\pm1$, as mentioned in the lines preceding
(\ref{gvc}). 

The question now is, if we exclude the prod\-uct-of-sur\-faces case (see
Remark~\ref{prdct}), can a function $\,Q\,$ of the form (\ref{qeq}) be
positive on an open interval while, at its endpoints, $\,Q=0\,$ and
the derivative $\,\dot Q\,$ has mutually opposite, nonzero values?

We will now prove Theorem~\ref{trhol} by answering this question in the
negative. Let us therefore suppose that such $\,Q\,$ exists.

Replacing the variable $\,\vt\in\bbR\smallsetminus\{\nh\gamma\}\,$ by 
$\,\theta=(\log|\vt-\gamma|)/2$, we get $\,Q\,$ equal to
$\,\delta\ve Ke^{2\theta}/2\,$ plus a linear combination of 
$\,e^\theta\nnh\cos\nnh\sqrt{7\,}\nnh\theta\,$ and
$\,e^\theta\nnh\sin\nnh\sqrt{7\,}\nnh\theta$,
where $\,\delta=\mathrm{sgn}\hs(\vt-\gamma)$. As
$\,d\theta/d\vt=\delta e^{-\nh2\theta}\nnh/2$, there exists 
an open interval of the variable $\,\theta\,$ on which $\,Q\,$ is positive 
while, at the endpoints, $\,Q=0\,$ and $\,e^{-\nh2\theta}\hn dQ/d\theta\,$
has mutually opposite, nonzero values. We have
$\,2e^{-\theta}\hn Q=\delta\ve Ke^\theta\nh
+p\hs\sin\nnh\sqrt{7\,}\nnh(\theta-q)\,$ for some constants $\,p,q$, and
$\,p\ne0$ (as $\,Q\,$ in neither constant, nor monotone).

Consequently, at some two distinct values of $\,\theta\,$ 
the function sending $\,\theta\in\bbR\,$ to 
$\,p\hs e^{-\theta}\nnh\sin\nnh\sqrt{7\,}\nnh(\theta-q)\,$ 
assumes the same value $\,-\delta\ve K$, with opposite nonzero values of the
derivative, and is greater than $\,-\delta\ve K\,$ between them. Note that we
are free to set $\,q=0$ by replacing the variable
$\,\theta\,$ with $\,\theta-q\,$ and $\,p\,$ with $\,p\hs e^{-q}\nh$.
Rescaling $\,K$, we may further assume that $\,|p|=1$. 
Thus, for $\,F:\bbR\to\bbR\,$ with 
$\,F(\theta)=e^{-\theta}\nnh\sin\nnh\sqrt{7\,}\nnh\theta$, there are 
\begin{itemize}
\item[$(*)$] two choices of $\,\theta\,$ at which $\,dF/d\theta\,$ has
opposite nonzero values, while $\,F\,$ assumes the same value, 
and is different from this last value between them.
\end{itemize}
Treating $\,F\,$ as a function of the new variable
$\,\alpha=\sqrt{7\,}\nnh\theta\,$ and setting
$\,c=\tan^{-\nh1}\hskip-3pt\sqrt{7\,}$, we get ($*$), with $\,\theta\,$
replaced by $\,\alpha$, for $\,F(\alpha)=e^{-\alpha\cot c}\sin\alpha$.
This contradicts Lemma~\ref{nexst}, thus completing the proof of
Theorem~\ref{trhol}.

\section{New examples of weakly Einstein K\"ah\-ler surfaces}\label{ne}
The construction mentioned in Theorem~\ref{weskr} uses any function 
$\,Q>0\,$ of the variable $\,\vt$ having the form (\ref{qeq}), that is, any
positive solution of (\ref{eul}), to define a weakly Einstein K\"ah\-ler
surface, which -- according to Remarks~\ref{divby} and ~\ref{notlh}
-- is neither Ein\-stein nor locally homogenous
unless $\,Q=\ve K(\vt-\gamma)/2$, for the constant Gauss\-i\-an curvature
$\,K\,$ of the base-sur\-face metric $\,h$.

For the reader's benefit we provide below a different, self-con\-tain\-ed
description of these examples, reflecting the fact that they have 
co\-ho\-mo\-ge\-ne\-i\-ty one (see Remark~\ref{chone} below), and
generalizing a construction in
\cite{grdq81}, rather than following the approach of \cite{plms03}.

In other
words, the goal of this section is to offer a more user-friend\-ly version
of the presentation given in Sect.\,\ref{ls}.

Let us fix nonzero real constants $\,\px,\qx\,$ and consider a
four-man\-i\-fold $\,M\,$ with vector fields $\,e_1\w,\dots,e_4\w$ 
trivializing $\,T\nh M\,$ and satisfying the Lie-brack\-et relations
\begin{equation}\label{eoe}
\begin{array}{l}
[e_1\w,e_i\w]=0\,\,\mathrm{\ for\ }\,i=2,3,4,\,\,\quad[e_2\w,e_4\w]
=2\px e_3\w,\\
{}[e_2\w,e_3\w]=\qx\hh e_4\w,\quad[e_3\w,e_4\w]=\qx\hh e_2\w.
\end{array}
\end{equation}
We define a metric $\,g\,$ and an al\-most-com\-plex structure $\,J\,$ on
$\,M\,$ by
\begin{equation}\label{geo}
g(e_1\w,e_1\w)=g(e_3\w,e_3\w)=\zeta\nh\eta,\quad
g(e_2\w,e_2\w)=g(e_4\w,e_4\w)=\zeta,\quad Je_1\w=e_3\w,\quad Je_2\w=e_4\w
\end{equation}
and $\,g(e\hn_i\w,e\nh_k\w)=0\,$ otherwise. 
Here $\,\zeta,\eta,\theta\,$ are functions of the real variable $\,\ut$, with
$\,\zeta,\eta$ assumed positive, and $\,\ut\,$ also stands for
a function $\,\ut:M\to\bbR\,$ such that
\begin{equation}\label{deo}
d_{e_1\w}\nnh \ut=2\zeta\nh\eta\hs\theta,\quad
d_{e_i\w}\nnh \ut=0\,\mathrm{\ for\ }\,i>1,\,\,\,
\mathrm{and\ so\ }\,\nabla\nnh\ut=2\theta\hs e_1\w,\quad
g(\nabla\nnh \ut,\nabla\nnh \ut)=4\zeta\nh\eta\hs\theta^2\nh.
\end{equation}
Such $\,\vt\,$ exists, locally, due to the obvious closedness of the
$\,1$-form $\,d\vt$, sending $\,e_1\w$ to $\,2\hh\zeta^2\nh\eta^2\hn\theta\,$
and the other three $\,e\hn_i\w$ to $\,0$.

Writing, this time, $\,(\,\,)'\nh=\hs d/d\ut$, we also assume that 
$\,\zeta'\nh\theta=-\px$, which turns out to guarantee that
$\,(M\nh,g,J)\,$ is a K\"ah\-ler manifold.

The geometric content of our discussion remains unchanged when $\,\ut\,$ is 
replaced by any function $\,\chi\,$ of the real variable $\,\ut$, via a
diffeomorphic change of the variable (that is, with $\,|\chi'|>0$). 
As functions on $\,M\nh$, our $\,\zeta\,$ and $\,\eta\,$ the remain the
same, while the role of $\,\theta\,$ is now played by 
$\,\theta\nh_{\mathrm{new}}\w=\chi'\nh\theta$. The function
$\,\zeta'\nh\theta\,$ and the condition $\,\zeta'\nh\theta=-\px\,$ 
remain unaffected. 
We use this freedom of modifying $\,\ut\,$ to require, without
loss of generality, that $\,\ut:M\to\bbR\smallsetminus\{\nh\gamma\}$, while
\begin{equation}\label{zet}
\theta\ne0\,\mathrm{\ \ be\ constant\ and\ \ }\,
\zeta\theta=\px(\gq-\ut)\,\mathrm{\ \,for\ some\ }\,\gq\in\bbR.
\end{equation}
Namely, as $\,\theta\ne0\,$ everywhere (due to the condition
$\,\zeta'\nh\theta=-\px$), we may choose $\,\chi\,$ above so as
to make $\,\theta\nh_{\mathrm{new}}\w$ constant.

Since $\,\zeta'\nh\theta=-\px$, 
the Le\-vi-Ci\-vi\-ta connection $\,\nabla\hs$ of $\,g\,$ is given by
\begin{equation}\label{neo}
\begin{array}{l}
\nabla\hskip-4pt_{e_1\w}\w\nnh\nh e_1\w=(\zeta\nh\eta)'\nh\theta e_1\w,\,\,\,
\nabla\hskip-4pt_{e_1\w}\w\nnh\nh e_2\w=\nabla\hskip-4pt_{e_2\w}\w\nnh\nh e_1\w
=-\px\eta e_2\w,\\
\nabla\hskip-4pt_{e_1\w}\w\nnh\nh e_3\w=\nabla\hskip-4pt_{e_3\w}\w\nnh\nh e_1\w
=(\zeta\nh\eta)'\nh\theta e_3\w,\,\,\,
\nabla\hskip-4pt_{e_1\w}\w\nnh\nh e_4\w=\nabla\hskip-4pt_{e_4\w}\w\nnh\nh e_1\w
=-\px\eta e_4\w,\\
\nabla\hskip-4pt_{e_2\w}\w\nnh\nh e_2\w=\px e_1\w,\,\,\,
\nabla\hskip-4pt_{e_2\w}\w\nnh\nh e_3\w=-\px\eta e_4\w,\,\,\,
\nabla\hskip-4pt_{e_2\w}\w\nnh\nh e_4\w=\px e_3\w,\\
\nabla\hskip-4pt_{e_3\w}\w\nnh\nh e_2\w=-(\px\eta+\qx)e_4\w,\,\,\,
\nabla\hskip-4pt_{e_3\w}\w\nnh\nh e_3\w=-(\zeta\nh\eta)'\nh\theta e_1\w,\\
\nabla\hskip-4pt_{e_3\w}\w\nnh\nh e_4\w=(\px\eta+\qx)e_2\w,\\
\nabla\hskip-4pt_{e_4\w}\w\nnh\nh e_2\w=-\px e_3\w,\,\,\, 
\nabla\hskip-4pt_{e_4\w}\w\nnh\nh e_3\w=\px\eta e_2\w,\,\,\,
\nabla\hskip-4pt_{e_4\w}\w\nnh\nh e_4\w=\px e_1\w.
\end{array}
\end{equation}
Also, as $\,\zeta'\nh\theta=-\px$, the only
pos\-si\-bly-non\-ze\-ro components 
of the curvature tensor $\,R$, the Ricci tensor $\,r\,$ 
and the metric $\,g\,$
are those algebraically related to
\begin{equation}\label{rot}
\begin{array}{l}
R_{1212}\w=R_{1234}\w=R_{1414}\w=R_{1423}\w=R_{2323}\w=R_{3434}\w
=\px\zeta^2\hn\eta\eta'\nh\theta,\\
R_{1313}\w=-2[(\zeta\nh\eta)'\nh\theta]'\nh\zeta^2\hn\eta^2\hn\theta,\quad
R_{1324}\w=2\px\zeta^2\hn\eta\eta'\nh\theta,\phantom{1^{1^1}}\\
R_{2424}\w=-\nh2\px(2\px\eta+\qx)\hs\zeta,\\
r_{11}\w=r_{33}\w=2[3\px\eta'\nh
-\zeta(\eta'\nh\theta)']\zeta\nh\eta\hs\theta,\phantom{1^{1^1}}\\
r_{22}\w=r_{44}\w=2\px\zeta\nh\eta'\nh\theta
-\nh2\px(2\px\eta+\qx),\phantom{1^{1^1}}\\
g_{11}\w=g_{33}\w=\zeta\nh\eta,\qquad g_{22}\w=g_{44}\w=\zeta.\phantom{1^{1^1}}
\end{array}
\end{equation}
For $\,\vz=\mathrm{Span}\hs(e_1\w,e_3\w)\,$ and
$\,\hz=\mathrm{Span}\hs(e_2\w,e_4\w)$, we thus get
\begin{equation}\label{rmg}
\ric=\mu g\,\mathrm{\ on\ }\,\mathcal{V}\nnh,\qquad\ric=\lambda g\,\mathrm{\
on\ }\,\mathcal{H},\qquad\ric(\mathcal{V}\nnh,\mathcal{H})=\{0\},
\end{equation}
where $\,\mu=r_{11}\w/g_{11}\w$ and $\,\lambda=r_{22}\w/g_{22}\w$.
Consequently, 
with $\,e_1\w\otimes e_1\w+e_3\w\otimes e_3\w$ equal to 
$\,g_{11}\w g\,$ on $\,\vz\hs$ and to $\,0\,$ on $\,\hz\nh$, we get,
as in (\ref{rre}),
\[
\ric\,=\,\lambda g\,+\,\frac{\mu-\lambda}{g_{11}\w}\hs(e_1\w\otimes e_1\w
+e_3\w\otimes e_3\w),
\]
provided that we identify $\,e_1\w\,$ with $\,g(e_1\w,\,\cdot\,)$, and
similarly for $\,e_3\w$. Therefore,
\begin{equation}\label{rri}
R\hh\ric\,=\,\lambda\ric\,
+\,\,\frac{\mu-\lambda}{g_{11}\w}\hs[R(e_1\w,\,\cdot\,,e_1\w,\,\cdot\,)
+R(e_3\w,\,\cdot\,,e_3\w, \,\cdot\,)].
\end{equation}
Since $\,R\hh\ric\,$ is $\,J$-in\-var\-i\-ant, and clearly
di\-ag\-o\-nal\-iz\-ed by our frame, $\,g\,$ is weakly Ein\-stein if and
only if $\,[R\hh\ric]_{11}\w/g_{11}\w=[R\hh\ric]_{22}\w/g_{22}\w$, that is,
by (\ref{rri}), with $\,R_{2323}\w=R_{1212}\w$ in (\ref{rot}),
\[
\lambda\mu\,+\,\frac{(\mu-\lambda)R_{1313}\w}{g_{11}^2}\,
=\,\lambda^2\hs+\,\frac{2(\mu-\lambda)R_{1212}\w}{g_{11}\w g_{22}\w}\hs.
\]
The Ein\-stein case being excluded, we can subtract the two sides and
-- according to Remark~\ref{divby} -- divide by $\,\lambda-\mu$, obtaining 
$\,\lambda+R_{1313}\w/g_{11}^2\nh=
2R_{1212}\w/[g_{11}\w g_{22}\w]$. As (\ref{rot}) implies that 
$\,R_{1313}\w/g_{11}^2\nh=-2[(\zeta\nh\eta)'\nh\theta]'\hn\theta\,$ and
$\,R_{1212}\w/[g_{11}\w g_{22}\w]=\px\eta'\nh\theta$, this reads
\[
\lambda\,-\,2[(\zeta\nh\eta)'\nh\theta]'\nh\theta\,
-\,2\hh\px\eta'\nh\theta\,=\,0
\]
and, multiplied by $\,g_{22}\w=\zeta$, it becomes
$\,r_{22}\w-2[(\zeta\nh\eta)'\nh\theta]'\nh\zeta\theta
-2\hh\px\hs\zeta\nh\eta'\nh\theta=0\,$ or, equivalently due to (\ref{rot}),
$\,2\px\zeta\nh\eta'\nh\theta-\nh2\px(2\px\eta+\qx)
-2[(\zeta\nh\eta)'\nh\theta]'\nh\zeta\theta
-2\hh\px\hs\zeta\nh\eta'\nh\theta=0$, that is
\[
\px(2\px\eta+\qx)\,
+\,[(\zeta\nh\eta)'\nh\theta]'\nh\zeta\theta\,=\,0.
\]
Since $\,\theta\,$ is constant and $\,\zeta\theta=\px(\gq-\ut)\,$ in
(\ref{zet}), multiplying by $\,4\zeta\theta^2$ we get
\begin{equation}\label{umq}
(\ut-\gq)^2Q\hs''\hs+\,2Q\,=\,4\qx\hh\theta(\ut-\gq),
\end{equation}
for $\,Q=4\zeta\nh\eta\hs\theta^2\nh$, 
which is precisely equation (\ref{eul}) with 
$\,K=4\ve\qx\hh\theta$, where $\,\ve=\pm\nh1\,$ is the signum of 
$\,\ut-\gq\,$ on our interval of the variable $\,\vt$, cf.\ (\ref{gvc}).
\begin{rem}\label{chone}The Lie-brack\-et relations (\ref{eoe}) define a
di\-rect sum Lie algebra $\,\bbR\hn\oplus\hh\mathfrak{g}$, for 
$\,\mathfrak{g}\,$ 
isomorphic to $\,\mathfrak{sl}\hh(2,\bbR)\,$ or
$\,\mathfrak{su}\hh(2)=\mathfrak{so}\hh(3)$, so that, locally,
$\,e_1\w,\dots,e_4\w$ are left-in\-var\-i\-ant 
vector fields on a di\-rect-prod\-uct Lie group $\,\bbR\nh\times G\,$ with 
$\,e_1\w$ tangent to $\,\bbR\,$ and $\,e_2\w,e_3\w,e_4\w$ to $\,G$.
Right-in\-var\-i\-ant vector fields on $\,G$, transplanted into
$\,\bbR\nh\times G$, commute with $\,e_1\w,\dots,e_4\w$ and 
are functional combinations of $\,e_2\w,e_3\w,e_4\w$ which, by
(\ref{deo}), makes them $\,g$-or\-thog\-o\-nal to the gradient 
$\,\nabla\nnh\ut$. Their flows thus preserve $\,\zeta,\eta\,$ and
$\,e_1\w,\dots,e_4\w$ and, consequently, the metric $\,g\,$ defined by
(\ref{geo}), so that $\,g\,$ has local co\-ho\-mo\-ge\-ne\-i\-ty one.
\end{rem}

\ \

{\it {Acknowledgments}}.
The research of YE was supported by Basic Science Research Program through
the National Research Foundation of Korea (NRF) funded by the Ministry of
Education (Grant number RS-2023-00244736). The research of SK was supported
by Basic Science Research Program through the National Research Foundation of
Korea (NRF) funded by the Ministry of Education (Grant number
RS-2023-00247409). The research of JP was supported by the National Research
Foundation of Korea (NRF) grant funded by the Korea government (MSIT)
(RS-2024-00334956). The authors also thank the anonymous referee for comments
and suggestions that greatly improved the exposition.



\begin{thebibliography}{99}

\bibitem{AK}
T.\,Arias-Marco and O.\,Kowalski, \textit{Classification of $4$-dimensional homogeneous weakly Ein\-stein manifolds}, Czechoslovak Math. J. \textbf{65}(140) (2015),
\hbox{21\hs--59}.

\bibitem{besse}
A.\,L.\,Besse, \textit{Ein\-stein manifolds}, Ergebnisse der Mathematik und ihrer Grenzgebiete (3), Springer-Verlag, Berlin, 1987.

\bibitem{BV}
E.\,Boeckx and L.\,Vanhecke, \textit{Unit tangent sphere bundles with
constant scalar curvature}, Czechoslovak Math.\,J. \textbf{51}(126)
(2001), 
\hbox{523\hs--544}.

\bibitem{bolton}J.\,Bolton, \textit{Transnormal systems}, Quart.\ J.\ Math.\ Oxford Ser. (2) \textbf{24} (1973), \hbox{385\hs--395}.

\bibitem{bourguignon}
J.\nh-P\nnh.\,Bourguignon, \textit{Les vari\'et\'es de dimension $4$ \`a signature non nulle dont la courbure est harmonique sont d'Ein\-stein}, 
Invent. Math. \textbf{63} (1981), \hbox{263\hs--286}.

\bibitem{calabi}E.\,Calabi, \textit{Extremal K\"ah\-ler metrics}, in: Yau, S.-T. (ed.), Seminar on Differential Geometry, 259--290. Annals of Math.\ 
Studies \textbf{102}, Princeton Univ.\ Press, Princeton, NJ, 1982. 

\bibitem{cartan}
E.\,Car\-tan, \textit{Sur une classe remarquable d'espaces de Riemann},
Bull. Soc. Math. France \textbf{l54} (1926), 
\hbox{214\hs--264}.

\bibitem{grdq81}
A.\,Derdzi\'nski, \textit{Exemples de m\'etriques de K\"ah\-ler et d'Ein\-stein autoduales sur le plan complexe}, G\'eo\-m\'e\-trie riemannienne
en dimension 4, (S\'eminaire Arthur Besse 1978/1979), Cedic/Fernand Na\-than, Paris (1981), pp. \hbox{334\hs--346}.

\bibitem{cm83} A.\,Derdzi\'nski, \textit{Self-dual Käh\-ler manifolds and Ein\-stein manifolds of dimension four}, Compos.\ Math. \textbf{49} (1983)
\hbox{405\hs--433}.

\bibitem{hdg00} A.\,Derdzinski, \textit{Ein\-stein metrics in dimension four},
Handbook of Differential Geometry, vol. I, pp. \hbox{419\hs--707}. North-Holland, Amsterdam (2000).

\bibitem{iumj12} A.\,Derdzinski, \textit{Kil\-ling potentials with geodesic gradients on K\"ah\-ler surfaces}, Indiana Univ.\ Math.\ J., \textbf{61}
(2012), 
\hbox{1643\hs--1666}.
  
\bibitem{plms03}A.\,Derdzinski and G.\,Maschler, \textit{Local 
classification of con\-for\-mal\-ly-Ein\-stein K\"ah\-ler metrics in higher 
dimensions}, Proc.\ London Math.\ Soc. \textbf{87} (2003), \hbox{779\hs--819}.

\bibitem{jram06} A. Derdzinski and G. Maschler, \textit{Special 
K\"ah\-ler-Ric\-ci potentials on compact K\"ah\-ler manifolds}, J.\ reine 
angew.\ Math. \textbf{593} (2006), 
\hbox{73\hs--116}.

\bibitem{EPS13}
Y.\,Euh, J.\,H.\,Park and K.\,Sekigawa, \textit{A curvature identity on a $4$-dimensional Riemannian manifold}, Results Math. \textbf{63} (2013), 
\hbox{107\hs--114}.

\bibitem{EPS14}
Y.\,Euh, J.\,H.\,Park and K.\,Sekigawa, \textit{A generalization of a $4$-dimensional Ein\-stein manifold}, Math. Slovaca \textbf{63} (2013), 
\hbox{595\hs--610}.

\bibitem{EPS15}
Y.\,Euh, J.\,H.\,Park and K.\,Sekigawa, \textit{Critical metrics for quadratic functionals
in the curvature on 4-dimensional manifolds}, Differen. Geom.
Appl. \textbf{29} (2011), \hbox{642\hs--646}.

\bibitem{GHMV}
E.\,Garc\'{i}a-R\'{i}o, A.\,Haji-Badali, R.\,Mari\~{n}o-Villar and M.\,E.\,V\'{a}zquez-Abal, \textit{Locally conformally flat weakly-Ein\-stein manifolds}, Arch. Math. (Basel) \textbf{111} (2018), 
\hbox{549\hs--559}.

\bibitem{GW}
A.\,Gray and T.\,J.\,Willmore, \textit{Mean-value theorems for Riemannian
manifolds},
Proc. Roy. Soc. Edinburgh Sect. A \textbf{92} (1982), 
\hbox{343\hs--364}. 

\bibitem{miyaoka}R.\,Miyaoka, \emph{Trans\-nor\-mal functions on a Riemannian 
manifold}. Differential Geom.\ Appl. \textbf{31} (2013), \hbox{130\hs--139}.

\bibitem{KNP} J.\,Kim, Y\nnh.\,Nikolayevsky and J.\,H.\,Park,
\emph{Weakly Ein\-stein Hypersurfaces in space forms}, preprint, available
from https:/\hskip-1.7pt/arxiv.org/abs/2409.12766v1.

\bibitem{tanno}
S.\,Tanno, \emph{$4$-dimensional conformally flat Kaehler manifolds},
 Tohoku Math.\,J. \textbf{24} (1972), \hbox{501\hs--504}.

\bibitem{wang}Q.\,M.\,Wang, \textit{Isoparametric functions on 
Riemannian manifolds}, I. Math.\ Ann. \textbf{277} (1987), 
\hbox{639\hs--646}.

\end{thebibliography}
\end{document}